\documentclass{article}

\usepackage[bottom]{footmisc}
\usepackage{geometry}
\usepackage{amssymb}
\usepackage{amsmath}
\usepackage{amsthm}
\usepackage{natbib}
\usepackage{theoremref}
\usepackage{todonotes}
\usepackage{tikz}
\usetikzlibrary{arrows}
\usepackage{caption}
\usepackage{subcaption}
\usepackage{float}
\usepackage{hyperref}

\usepackage{doi}

\usepackage{times}
\usepackage{bm}
\usepackage{natbib}
\usepackage{tikz-cd}
\usepackage{url}
\usepackage{mathtools}
\usepackage[capitalise,noabbrev]{cleveref}

\usepackage{authblk}

\graphicspath{{./art/}}

\usepackage{adjustbox}
\usepackage{makebox}

\usepackage{tikz-cd}
\usepackage{mathdots}

\usepackage[ruled,vlined]{algorithm2e}
\SetKwProg{init}{Initialization}{}{}
\usepackage{setspace}

\newtheorem{definition}{Definition}
\newtheorem{theorem}{Theorem}
\newtheorem{lemma}{Lemma}

\newtheorem{remark}{Remark}

\newcommand\numberthis{\addtocounter{equation}{1}\tag{\theequation}}
\DeclareMathOperator*{\argmax}{argmax}
\DeclareMathOperator*{\argmin}{argmin}
\SetKwProg{init}{Initialization}{}{}
\def\CCT{{Cauchy Combination Test}}
\def\CBAM{{Compound Bonferroni Arithmetic Mean}}
\def\HMP{{Harmonic Mean $p$-value}}
\def\FCT{{Fisher Combination Test}}
\def\TPM{{Truncated Product Method}}
\def\rTPM{{Rank Truncated Product Method}}
\def\CTPs{{Closed Testing Procedures}}
\def\CTP{{Closed Testing Procedure}}
\def\prop{{\mathbb{P}}}  %

\newcommand{\changed}[1]{{#1}}
\usepackage{enumitem}
\newlist{todolist}{itemize}{2}
\setlist[todolist]{label=$\square$}
\usepackage{pifont}

\begin{document}

\title{%
Too Many, Too Improbable: testing joint hypotheses and closed testing shortcuts
}
\author{Phillip B. Mogensen and Bo Markussen}
\affil[1]{Department of Mathematical Sciences, University of Copenhagen\\ Universitetsparken 5, Copenhagen, Denmark}
\date{}

\maketitle

\begin{abstract}
Hypothesis testing is a key part of empirical science and multiple testing as well as the combination of evidence from several tests are continued areas of research. In this article we consider the problem of combining the results of multiple hypothesis tests to i) test global hypotheses and ii) make marginal inference while controlling the $k$-FWER. 
We propose a new family of combination tests for joint hypotheses, called the ‘Too Many, Too Improbable’ (TMTI) statistics, which we show through simulation to have higher power than other combination tests against many alternatives. Furthermore, we prove that a large family of combination tests -- which includes the one we propose but also other combination tests -- admits a quadratic shortcut when used in a  {\CTP}, which controls the FWER strongly. We develop an algorithm that is linear in the number of hypotheses for obtaining confidence sets for the number of false hypotheses among a collection of hypotheses and an algorithm that is cubic in the number of hypotheses for controlling the $k$-FWER for any $k$ greater than one.
\end{abstract}

\section{Introduction}

The problem of combining $p$-values from tests of a family of hypotheses $\{H_i\}_{i \in \mathcal I}$ indexed by a finite set $\mathcal{I}$ has long been a field of study and remains an active area of research today. In 1925, Fisher proposed a method of combining independent $p$-values by observing that minus two times the sum of $(\log p_i)_{i \in \mathcal I}$ follows a $\chi_{2\vert\mathcal I\vert}^2$-distribution \citep{fisher1992statistical}. Fisher's combination test is asymptotically Bahadur-optimal among the class of all combination tests \citep{littell1973asymptotic}. Still, the {\FCT}  can potentially be outperformed by other combination tests for any given finite sample. In 1973, Brown devised an extension of the {\FCT}  or when the underlying test statistics are jointly Gaussian with a known covariance matrix and the hypotheses are one-tailed \citep{brown1975400}. \citet{kost2002combining} further relaxed the assumptions on the dependence structure by deriving an approximation of the distribution of the {\FCT}  when the underlying tests statistics are jointly $T$-distributed with a common denominator. In recent years the {\CCT}  \citep{acat_original} and the {\HMP} \citep{wilson2019harmonic} have been proposed and \citet{vovk2020combining} derive a large family of combination tests by using the Kolmogorov generalized $f$-mean. 

These combination-based methods for testing the global null hypothesis $H_0 \coloneqq \bigcap_{i \in \mathcal I} H_i$ follow the overall recipe of finding a mapping $f: [0, 1]^{\vert\mathcal I\vert} \to [0,1]$ of $\boldsymbol{p} \coloneqq (p_i)_{i \in \mathcal I}$ such that $f(\boldsymbol{p})$ is again a $p$-value under $H_0$. That is, such  that for any choice of $\alpha \in [0,1]$ it holds that $\prop(f(\boldsymbol{p}) \leq \alpha) \leq \alpha$ when $H_0$ is true. One particular way of obtaining this property is to choose \textit{any} function, say $f_1$, that maps the hypercube $[0,1]^{\vert\mathcal I\vert}$ to any subset of the real line and then subsequently transform the resulting random variable by its cumulative distribution function (CDF), say $f_2$. The composite mapping $f_2 \circ f_1$ is then a valid combination test. The {\FCT}  is an example of this; first, we map $\boldsymbol{p}$ onto the positive real line by the mapping $\boldsymbol{p} \mapsto -2 \sum_{i \in \mathcal I} \log p_i$, which is then transformed back onto the unit interval using the CDF of a $\chi_{2\vert\mathcal I\vert}^2$-distribution. Another simple way of constructing valid combination tests is to use the minimal $p$-value from any procedure that controls the family-wise error rate (FWER). For example, we may use the minimal $p$-value of the Bonferroni corrected $p$-values, corresponding to the mapping $\boldsymbol{p} \mapsto \min\left({\vert\mathcal I\vert\cdot \min(\boldsymbol{p}), 1}\right)$. 

\changed{In this paper, we introduce a family of combination tests -- the ‘Too Many, Too Improbable’ (TMTI) tests -- that strongly controls the Type I error at level $\alpha$, for any choice of $\alpha \in (0, 1)$.}
In brief, these statistics are obtained by ordering the observed $p$-values, transforming them by the CDFs of beta distributions and returning a local minimum. The $p$-value is then the local minimum transformed by its CDF. We derive analytical expressions for the null CDFs of the TMTI test statistics under an assumption of independence and show through simulation that the TMTI tests can have higher power than other common combination tests under many alternatives. Additionally, we give an $\mathcal{O}(m^2)$ shortcut for carrying out a full {\CTP} for all elementary hypotheses for a large family of test statistics, obtained by considering test statistics of the form $Z = h(F_{(1)}(p_{(1)}), \dots, F_{(m)}(p_{(m)}))$ under mild assumptions on the functions $F_{(1)}, \dots, F_{(m)}$ and $h$.
Using prior work by \citet{goeman2011multiple}, we show how these shortcuts can be used to obtain $k$-FWER control for elementary hypotheses as well as construct confidence sets for the number of false hypotheses in a rejection set. Finally, we discuss how mixing different local tests across a  {\CTP}  can be used to increase power.

\section{The ‘Too Many, Too Improbable’ family of test statistics}\label{sec:setup}

\subsection{Notation and setup}
Let $\mathcal I \coloneqq \{1, \dots, m\}$ be an index set with cardinality $m$ and let $\{H_i\}_{i \in \mathcal{I}}$ be hypotheses.
Let $(P_i)_{i \in \mathcal{I}}$ be random variables on probability spaces $(\Omega_i, \mathbb B_i, \mathbb P_i)_{i \in \mathcal I}$ with $\Omega_i \subseteq [0,1]$. In most situations we will have $\Omega_i = [0, 1]$ and have $\mathbb B_i$ be the Borel sigma-algebra, although this need not be the case. We denote by $p_i$ an outcome of $P_i$ and call $p_i$ the $p$-value for the test of $H_i$. For a given subset of indices, $\mathcal J \subseteq \mathcal I$, we consider the task of testing the joint hypothesis
$
  H_{\mathcal J} \coloneqq \bigcap_{j \in \mathcal J} H_j
$
by using the marginal $p$-values, $\boldsymbol{P}^{\mathcal J} \coloneqq (P_j)_{j \in \mathcal J}$. The set $\mathcal J$ can be chosen freely according to what kind of hypothesis one wishes to test. If we choose $\mathcal{J}$ with $\vert\mathcal{J}\vert=1$, no adjustment needs to be made, as we are simply testing a marginal hypothesis, for which we already have a $p$-value. If we choose $\mathcal{J} = \mathcal{I}$, we are considering the global null hypothesis of $\mathcal{I}$. Anything in between those two extremes corresponds to testing a particular joint hypothesis. E.g., if $(p_i)_{i \in \mathcal{I}}$ are the $p$-values output from a genome-wide association study, then $\mathcal{J}$ could correspond to a particular region, which is of special interest, e.g., a gene or chromosome.

In order to test $H_\mathcal{J}$, we construct a test statistic, denoted by $Z$, with corresponding $p$-value $P$ that satisfies
	\begin{equation} \label{eq:level}
  		H_{\mathcal J} \text{ true} \implies \forall \alpha \in [0, 1]:~ \prop(P \leq \alpha) \leq \alpha.
	\end{equation}
The above statement is called Type I error control and means that whenever the joint hypothesis $H_{\mathcal J}$ is true, the probability that we reject $H_{\mathcal J}$ at level $\alpha$ is at most $\alpha$.

\subsection{Definition of the TMTI statistics}
Let $P_{(1)}^{\mathcal J} \leq \dots, P_{(\vert\mathcal J\vert)}^{\mathcal J}$ denote an ordering of $\boldsymbol{P}^{\mathcal J}$. This ordering is possibly not unique. Let $\beta(a, b)(x)$ denote the CDF of the $\beta(a,b)$-distribution.When the shape and scale parameters are integers, we can write
\begin{equation}\label{eq:pbeta_alt}
		\forall i, m \in \mathbb N: \hspace{1cm}
	\beta(i, m + 1 - i)(x) = \sum_{k = i}^{m} \binom{m}{k} x^k ( 1-x)^{m - k}.
\end{equation}%
We construct the collection $\boldsymbol{Y}^{\mathcal J} \coloneqq (Y^{\mathcal{J}}_k)_{k = 0}^{\vert\mathcal J\vert + 1}$ of random variables by $Y_0^{\mathcal J} \coloneqq 2$, $Y_{\vert\mathcal J\vert + 1}^{\mathcal J} \coloneqq 2$, and
\begin{displaymath}\label{eq:defineY}
	  \forall k \in \{1, \dots, \vert\mathcal J\vert\}: \hspace{1cm}
  			Y_k^{\mathcal J} \coloneqq \beta(k, \vert\mathcal J\vert + 1 - k) (P_{(k)}^{\mathcal J}).
\end{displaymath}
If all variables in $\boldsymbol{P}^{\mathcal{J}}$ are independent and exactly uniform, then each $Y^{\mathcal{J}}_k$ is uniformly distributed on $[0, 1]$ for $k=1,\dotsc,\vert\mathcal{J}\vert$, as it is well known that the order statistics of i.i.d.\ $U(0, 1)$ variables are $\beta$-distributed. 

Let $c \leq \vert \mathcal{J} \vert$ be an integer. 
We then consider the first $Y^{\mathcal{J}}_k$ among the first $c$ variables that is strictly smaller than the following $n \in \mathbb N \cup \{\infty\}$, i.e.,
\begin{displaymath}
  L_{n, c} \coloneqq \min \{l \in \{1, \dots, c\}: Y_l^{\mathcal J} < Y_{k}^{\mathcal J} ~~ \text{for all} ~~ k=l+1, \dots, \min(l + n, \vert\mathcal J\vert + 1) \}.
\end{displaymath}
If $Y_1^{\mathcal{J}} \geq \dots \geq Y_c^{\mathcal{J}}$ we set $L_{n,c} = c$.
We think of $L_{n, c}$ as the index of a kind of local minimum of $Y_1^{\mathcal{J}} , \dots , Y_c^{\mathcal{J}}$, in the sense that $Y_{L_{n, c}}^{\mathcal{J}}$ is always a local minimum, but it further needs to satisfy that it is smaller than the following $n$ terms. In particular, 
$Y_{L_{1, c}}^\mathcal{J}$ is the first local minimum of $Y_1^{\mathcal{J}} , \dots , Y_c^{\mathcal{J}}$ and 
$Y_{L_{\infty, c}}^{\mathcal J}$ is the global minimum of $Y_1^{\mathcal{J}} , \dots , Y_c^{\mathcal{J}}$. 
The construction of $Y^{\mathcal{J}}_0$ and $Y^{\mathcal{J}}_{\vert\mathcal{J}\vert + 1}$ is a technical one, meant only to ensure the existence of $L^{\mathcal{J}}_{n, c}$.
To ease the notational burden, we omit the subscripted $n$ and $c$ and the superscripted $\mathcal{J}$ when the particular choices of $n$, $c$ and $\mathcal J$ are not of importance or unambiguous from context. 

\begin{definition}
Let $n \in \mathbb{N} \cup \{\infty \}$ and let $c \leq \vert\mathcal{J}\vert$ be an integer. The ‘Too Many, Too Improbable’ test statistic is then defined as
\[
Z_{n, c}^{\mathcal J} \coloneqq Y_{L_{n, c}}^{\mathcal J}.
\]
Small values of $Z_{n, c}^{\mathcal J}$ are critical and the $p$-value for the test of $H_{\mathcal{J}}$ is obtained by evaluating the test statistic in its CDF under $H_{\mathcal{J}}$. We denote by $\gamma_{n, c}^{\mathcal{J}}(x)$ the CDF of $Z_{n, c}^{\mathcal J}$ under $H_{\mathcal{J}}$.
\end{definition}

\changed{%
Generally, we will only consider cases in which $n = 1$ or $n = \infty$, as these are the most natural choices of $n$. However, the setup allows for other choices of $n$. 
Choosing $1 < n < \infty$ can potentially increase the power of the procedure in cases where signals are fairly sparse, but sufficiently weak that the first local minimum falls ‘too early’ by chance. However, we do not investigate this further, but simply remark that it is possible to choose $n$ different from what we consider in the remainder of this paper.
}

\changed{%
\begin{remark}
    Testing the joint hypothesis $H_\mathcal{J}$ using any TMTI test satisfies the statement in \eqref{eq:level} by the probability integral transform. That is, the TMTI test controls the Type I error.
\end{remark}
}

\begin{remark}
	Whenever $\vert\mathcal J\vert = 1$, say $\mathcal J = \{j\}$, the TMTI transform is simply the identity transform, i.e., $\gamma(Z) = P_j$.
\end{remark}

\begin{remark}
	If the variables in $\boldsymbol{P}^\mathcal{I}$ are exchangeable, i.e., if any two subsets of equal size have the same joint distribution, it follows, that for any two sets $\mathcal J_1, \mathcal J_2 \subset \mathcal I$ with $\vert\mathcal J_1\vert = \vert\mathcal J_2\vert$, we have $\gamma^{\mathcal{J}_1} = \gamma^{\mathcal{J}_2}$. Thus, for exchangeable $p$-values, the CDF of the TMTI statistic depends only on the choice of $\mathcal J$ through its cardinality.
\end{remark}

\subsection{Truncation procedures and the TMTI}
The {\TPM}  of \citet{zaykin2002truncated} and the {\rTPM}  of \citet{dudbridge2003rank} are two notable variants of the {\FCT}, that also test the global null hypothesis but against different alternatives. 

The {\TPM}  is a combination test that uses only the $p$-values that are smaller than some predefined threshold $\tau \in (0, 1)$. The alternative hypothesis is therefore, that there is at least one false hypothesis among those hypotheses that gave rise to $p$-values below $\tau$. The {\rTPM}  is also a combination test, but this uses only the smallest $K$ $p$-values, for some predefined $K < \vert\mathcal{J}\vert$. Thus, the alternative hypothesis is, that there is at least one false hypothesis among those, that gave rise to the $K$ smallest $p$-values.

The TMTI family of test statistics includes similar procedures. For any $c < \vert \mathcal{J} \vert$, the alternative hypothesis is that there is at least one false hypothesis among those that gave rise to the $c$ smallest $p$-values. Thus, setting $c = K$ for some integer $K < \vert\mathcal{J}\vert$, the TMTI procedure uses only the first $K$ $p$-values in the construction of the test statistic and therefore tests the joint hypothesis $H_{\mathcal{J}}$ against the same alternative as the {\rTPM}. We call this procedure the rank truncated TMTI.%

By setting $c = \max(\{j \in \{1, \dots, \vert\mathcal{J}\vert\}: p_{(j)} \leq \tau \} \cup \{1\})\eqqcolon \bar{\tau}$, for some value $\tau \in (0, 1)$, the TMTI procedure uses only the $p$-values that are marginally significant at level $\tau$, and thus tests against the same alternative as the {\TPM}. We call this procedure the truncated TMTI.
In the event that no $p$-values are smaller than $\tau$, $c$ becomes $1$ and uses instead the smallest of the available $p$-values. 

We write TMTI$_n$ to denote the TMTI statistic $Z_{n, c}^{\mathcal{J}}$ with $c = \vert\mathcal{J}\vert$, tTMTI$_n$ to denote the truncated TMTI statistic and rtTMTI$_n$ to denote the rank truncated TMTI statistic.

\changed{%
There are two potential advantages to using a truncated procedure (i.e., $c < \vert \mathcal{J} \vert$) over a non-truncated procedure. First, for large $m$ (say, $m \geq 10^6$), it is non-trivial to compute the TMTI$_\infty$-statistic, because its computation involves sorting $m$ different $p$-values and computing $m$ different $\beta$-transformations. Using a truncation procedure instead reduces the computational cost, because fewer $p$-values need to be considered. Thus, only a partial sorting is required and fewer $\beta$-transformations need to be computed. Second, as we outline below, using a truncation procedure can potentially have higher power than its non-truncated version.
}

\begin{figure}[ht]
    \centering
  \includegraphics[height = 5cm, width = 10cm]{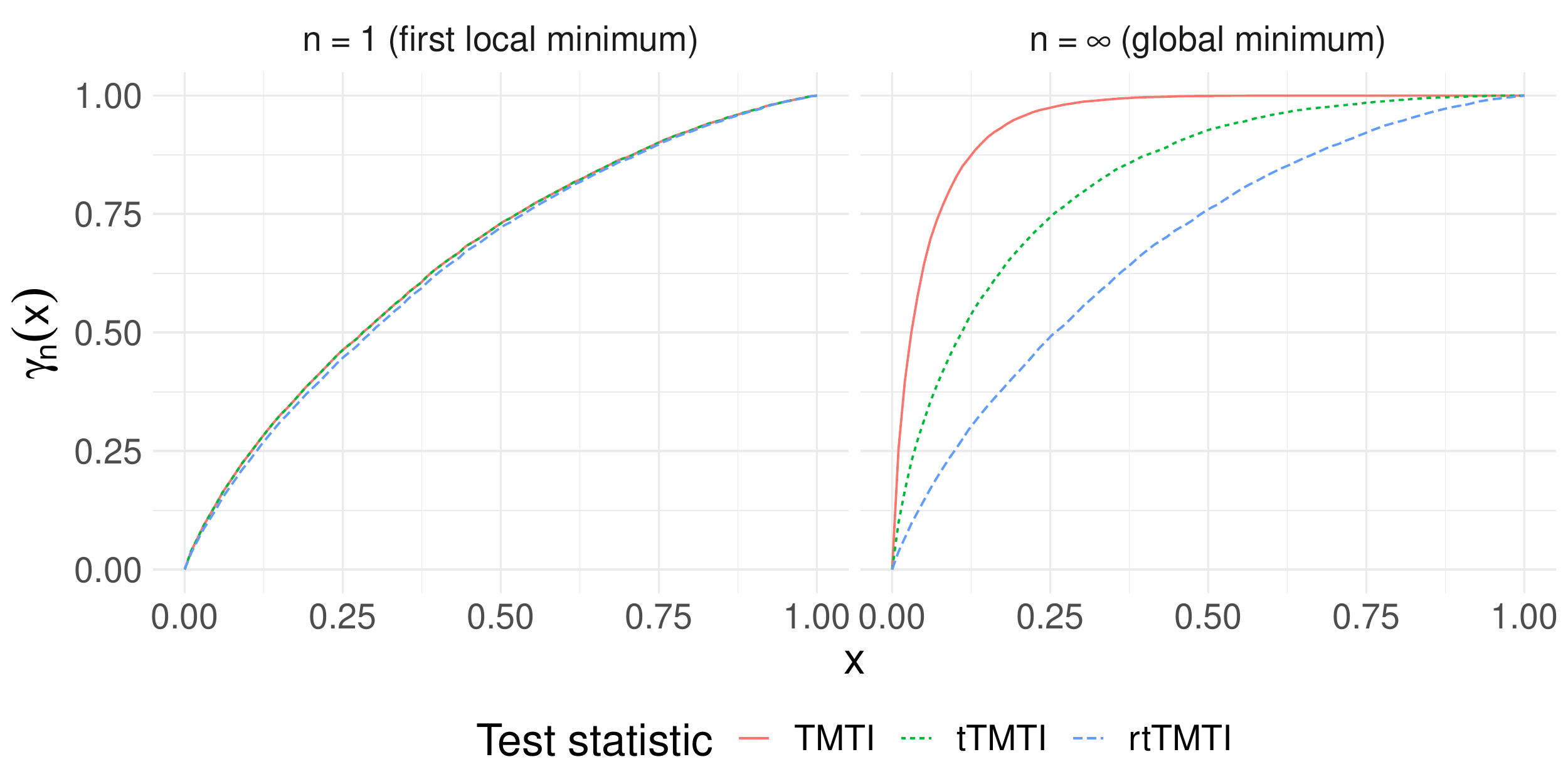}
  \caption{Comparison of $\gamma_{n}^\mathcal{I}$, $\gamma_{n, \tau}^\mathcal{I}$ and $\gamma_{n, K}^\mathcal{I}$ for $m = 10^5$, $\tau = 0.01$ and $K = 5$ in the case of independent and exactly uniform $p$-values. The solid red lines are TMTI, the dotted green and blue lines are tTMTI and rtTMTI, respectively.}
  \label{fig:gamma_comparison}
\end{figure}

\begin{lemma}\label{lemma:dominated}
	Let $n \in \mathbb{N} \cup \{\infty \}$ and $x \in (0, 1)$. Let $\mathcal{I}$ be an index set with cardinality $m$. It follows that 
	\[
	\forall c < m: \hspace{1cm}
	\gamma_{n, c}(x) < \gamma_{n,m} (x)
	\]
\end{lemma}
When using $n=1$, i.e., considering the first local minimum, and when using moderate values of $\tau$ and $K$, it is likely that $\gamma_{1,\tau}$, $\gamma_{1,K}$ and $\gamma_{1}$ are going to be nearly identical, as the first local minimum is likely to lie early in the sequence $\boldsymbol{Y}$. This implies that the $p$-values of the tTMTI$_1$, the rtTMTI$_1$ and TMTI$_1$ tests are nearly identical. Thus, the TMTI$_1$ by itself can be thought of as a truncation method. However, if using the global minimum, we expect that there can be a large difference between the methods, as the global minimum is likely to lie further along the sequence $\boldsymbol{Y}$. Thus, applying either the tTMTI$_\infty$ with a low $\tau$ or rtTMTI$_\infty$ with a low $K$ is going to be roughly equivalent to applying the TMTI$_1$. These properties are demonstrated in \cref{fig:gamma_comparison} for the case of independent and exactly uniform $p$-values.

From \cref{fig:gamma_comparison} we conclude, that if the global null hypothesis is indeed false, and if the global minimum of the sequence $\boldsymbol{Y}$ happens to fall within the first $K$ or $\bar{\tau}$ indices of $\boldsymbol{Y}$, there is a potential for a large power gain by applying either the tTMTI$_\infty$ or rtTMTI$_\infty$ over the standard TMTI$_\infty$. This is because the procedures will all be considering the same test statistic, $Z$, but they will evaluate it under different $\gamma$ functions, thereby yielding different $p$-values. By \cref{lemma:dominated}, the $\gamma$ functions from the tTMTI$_\infty$ and rtTMTI$_\infty$ procedures are dominated by the $\gamma$ function from the TMTI$_\infty$, implying that $p$-value resulting from applying the truncation procedures will be lower than those of the standard procedure, giving rise to higher power.

\section{Computation of $\gamma$}\label{sec:computation}

\subsection{An analytical expression of the CDF of TMTI$_\infty$ in the i.i.d. case} \label{sec:analytic_computation}
In the case where the $p$-values are independent under the null hypothesis, we can derive an analytical expression of $\gamma_{\infty, c}$.

\begin{theorem}\label{theorem:gamma}
	Let $P_1, \dots, P_m$ be i.i.d.\ uniformly distributed on $[0, 1]$. 
	For every $i \in \{1, \dots, m\}$, let 
	$x_i$ 
	be the $x$-quantile of the $\beta(i, m + 1 - i)$ distribution and define the polynomial
    \begin{equation*}%
        Q_{i}(x; \boldsymbol{a}) \coloneqq \sum_{j=1}^{i} \frac{a_j}{(i + 1 - j)!} x^{i + 1 - j}.
    \end{equation*}
    Define 
    $\Bar{Q}_1 \coloneqq x_1$ 
    and define recursively
    \begin{align*}
    \forall i \in \{2,\dots, m\}: \hspace{1cm}
    \Bar{Q}_{i, c} &\coloneqq Q_{i}(
    x_{\min(i, c)}
    ; (1, -\Bar{Q}_{1, c}, \dots, -\Bar{Q}_{i-1, c} )).
    \end{align*}
    If $c \leq m$ is a fixed integer, then
	\begin{equation}\label{eq:gamma_analytical}
    \gamma_{\infty, c}(x) = 
    x_c^m 
    + \sum_{i=1}^{m-1} \frac{m!}{(m - i)!} \bar{Q}_{i,c}  
    (1 - x_c^{m-i}).
	\end{equation}
Furthermore, let $\tau \in (0, 1)$ and define $\tilde{Q}_1 \coloneqq \tau$ and
	\[
	\forall i \in \{2, \dots, m\}: \hspace{1cm}
	\tilde{Q}_i \coloneqq Q_i(\tau; (1, -\bar{Q}_{1, m}, \dots, \bar{Q}_{i - 1, m})).
	\]
If $c$ is a random variable given by $c = \max(\{i \in \{1, \dots, m\}: P_{(i)} < \tau\}\cup \{1\})$, then
	\begin{equation}\label{eq:gamma_analytical_2}
	\begin{split}
		\gamma_{\infty, c}(x) &= 
		(1 - \tau)^m \frac{x - \beta(1, m)(\tau)}{1 - \beta(1, m)(\tau)} I(x_1 > \tau) \\
		&+ 
		\sum_{i = 1}^{m} \left[
			\binom{m}{i} \tau^i (1 - \tau)^{m - i} \left\{
				1 - \frac{i!}{\tau^i} \left(
					\tilde{Q}_i
					-
					\bar{Q}_{i, m}
				\right) 
				I(x_i \leq \tau)
			\right\}
		\right]
	\end{split}
	\end{equation}
\end{theorem}
The above can be readily implemented by recursively computing the $\bar{Q}_{i, c}$ and $\tilde{Q}_i$ terms, e.g.,in a for-loop. 

In the special case of $c = 1$, we have $\gamma_{n, 1}(x) = x$ by construction, regardless of $m$ and $n$. In this setting, the TMTI procedure is then a minimum-$p$ method. This has the advantage, that the procedure can be applied directly in high-dimensional settings, if the assumption of independence holds.
\changed{%
In particular, it is easy to show, that when $c = 1$ the critical value of the TMTI test is $1 - (1 - \alpha)^{1/m}$, and it is thus equivalent to using the Šidák correction \citep{vsidak1967rectangular} for testing the global null hypothesis. 
}

\subsection{A bootstrap scheme for the CDF of TMTI$_n$ in the i.i.d. case}\label{sec:bootstrap}
Although it is easy to implement Equations \eqref{eq:gamma_analytical} and \eqref{eq:gamma_analytical_2}, numerical difficulties may arise when $m$ is large, say $m > 100$, due to the presence of the factorials $1!, \dots, m!$ in the computations. Essentially, the $\bar{Q}_{i, c}$ terms are all very small, because they include the reciprocals of factorials, but they are scaled up by another factorial. Although this is well-defined, numerical instabilities will occur in implementations in standard double-precision arithmetic. For larger $m$, one can perform the calculations in arbitrary precision, but the added computational cost of doing so can be enormous. Instead, a simple bootstrapping scheme can be employed by; i) drawing $m$ values independently from a $U(0, 1)$ distribution; ii) transforming the values from step i as described in Equation \ref{eq:defineY} and saving the desired TMTI statistic as $Z_{b}$, where $b$ indexes the iteration; iii) repeating steps i and ii sufficiently many times, say $B$, and; iv) using $\hat{\gamma} (x) \coloneqq \frac{1}{B}\sum_{b=1}^B I(Z_b \leq x)$ as an approximation of $\gamma$. 
This bootstrap scheme can be applied regardless of the choice of $n$ and $c$.

\subsection{The non-independent case}
\changed{%
The derivation of the CDF, $\gamma$, in \cref{sec:analytic_computation} relied on the assumption that all $p$-values are independent. \citet{chen2020trade} argue that combining methods that are Valid under Arbitrary Dependence (VAD) structures have lower power than combining methods that are Valid under Independent (VI) $p$-values, if the underlying, $p$-values are in fact independent. However, if the underlying $p$-values are not independent, VI methods may fail to hold level, whereas their VAD counterparts will hold level for any dependence structure. Thus, the choice of combining method should depend on the scientific question of interest. For instance, in genome-wide association studies (GWAS), it is unreasonable to assume independence, as base-pairs are likely to be locally dependent \citep{dudbridge2003rank}.\footnote{It is often possible to filter the $p$-values in a manner such that the remaining $p$-values are likely to be independent, e.g., using a distance-based filtering.}

In some cases, however, it is possible to apply the methods described in \cref{sec:analytic_computation,sec:bootstrap}, even if the $p$-values are not independent. For instance, if the underlying tests, $Z$, are jointly Gaussian with a known covariance matrix, $\Sigma$, these can be decorrelated by performing an eigendecomposition, $\Sigma = Q \Lambda Q^{T}$, and then constructing $\tilde{Z} \coloneqq (Q \Lambda^{-1/2} Q^T)^T Z$. Then, the components of $\tilde{Z}$ are jointly independent (see, e.g., \citet{kessy2018optimal}) and the methods described in \cref{sec:analytic_computation,sec:bootstrap} can be directly applied to the $p$-values obtained from the test statistics $\tilde{Z}$. This remains true for any rotation, $\bar{Z} \coloneqq R (Q \Lambda^{-1/2} Q^T)^T Z$, where $R$ is an orthogonal matrix.

If the $p$-values are not independent, and if the decorrelation procedure described above is not appropriate, one can still try to apply the TMTI directly. However, level of the test (i.e., \cref{eq:level}) is no longer guaranteed, and thus there is a chance that the Type I error is increased. How much the Type I error increases depends entirely on the dependence structure of the $p$-values. In \ref{appendix:dependency_simulation}, we investigate the level of the TMTI tests under three different types of dependencies: autoregressive $p$-values (i.e., $\operatorname{cor}(P_i, P_j) = \rho^{\vert i - j\vert}$), equicorrelated $p$-values (i.e., $\operatorname{cor}(P_i, P_j) = \rho$, for all $i,j$), and block-diagonally correlated $p$-values (i.e., $\operatorname{cor}(\boldsymbol{P})$ has a block-diagonal structure, where all off-diagonal entries are $\rho$ if in the same block and $0$ else). We note, however, that rtTMTI$_n$ seems to either have the correct level or be conservative, no matter the dependence structure.
Overall, the TMTI tests hold level only under weak autoregressive and block-diagonal dependency structures, and fails to hold level for stronger dependency structures and equicorrelated $p$-values (see \ref{appendix:dependency_simulation} and \cref{fig:dependencyplot,fig:dependencyplot_zoom} for a full account of the results). Thus, the TMTI tests can potentially still be applied in settings with weak dependence, but it is not appropriate in settings with strong dependencies.

Finally, one can apply any VI combination method under arbitrary dependence, if one is able to sample from the joint distribution of $\boldsymbol{P}$ under the global null. This can, for instance, be done if one assumes an underlying parametric model or by employing a resampling bootstrap procedure. In \ref{sec:power}, we give an example of how this can be done in a case where the marginal hypotheses of interest are $T$-tests for the parameters in a linear model being zero.

}

\changed{%
\section{Power of the TMTI -- a simulation study}\label{sec:newpower}
In this section, we show by means of simulation, that many of the TMTI tests have high power against a wide range of alternative hypotheses. In particular, we find that the TMTI$_\infty$ and tTMTI$_\infty$ tests have high power both in cases where signals are sparse but strong and in cases where signals are dense but weak.

We consider $m = 10^5$ independent tests, of which $N_{\operatorname{false}} \in \{10^0, \cdots, 10^4\}$ are false. 
In order to investigate situations in which the $p$-values from true hypotheses are conservative, we generate these as $p_{\operatorname{true}} \coloneqq U^\delta$, where $U \sim U(0, 1)$ and $\delta \in [0, 1]$. 
When $\delta = 1$, this corresponds to the true $p$-values being exactly uniform, and when $\delta \in (0, 1)$, it corresponds to the true $p$-values being strictly conservative. The degree of conservatism increases as $\delta$ decreases and the extreme case in which $\delta = 0$ corresponds to the degenerate case where all true $p$-values are equal to one, meaning that no hypothesis can ever be rejected, no matter the significance level.
Situations in which the $p$-values are strictly conservative occur in many places. For instance, in a GWAS with dichotomous traits, the $p$-values will be conservative \citep{wu2011rare}. Conservative $p$-values also occur in Invariant Causal Prediction, where a $p$-value for invariance is obtained as the minimum of Bonferroni-corrected $p$-values from multiple environments \citep{peters2016causal}.
We generate $p$-values for the false hypotheses by independently sampling $Z$-scores, $Z_1, \dots, Z_{n_{\operatorname{false}}} \stackrel{i.i.d.}{\sim} N(\mu_{\operatorname{false}}, 1)$, for different values of $\mu_{\operatorname{false}}$ and then letting $p_{i, \operatorname{false}} \coloneqq 2 \times (1 - \Phi(\vert Z_i \vert))$, where $\Phi$ is the CDF of a $N(0, 1)$-distribution. The values of $\mu_{\operatorname{false}}$ are chosen equidistantly between the two values, which satisfy that a Bonferroni test has either 5\% or 99\% power to reject the global null hypothesis in a setting with no conservatism.

For comparison, we include the {\FCT} (which is known to lose power in the presence of conservative $p$-values \citep{zaykin2002truncated}) and its truncated versions, and the {\CCT}  and {\HMP} (which are known to have high power in settings with sparse, strong signals \citep{acat_original,wilson2019harmonic}). In all settings, we employ a significance level of $\alpha = 0.05$, and for the truncation procedures, we use $\tau = 0.05$ and $K = 10$. We include TMTI$_\infty$ and both truncation variants, as well as TMTI$_1$. For TMTI$_1$, we do not include any truncation variants, as we expect these to be roughly equal to the non-truncated version (per \cref{fig:gamma_comparison}).
\begin{figure}[t]
    \centering
    \includegraphics[width = .9\linewidth]{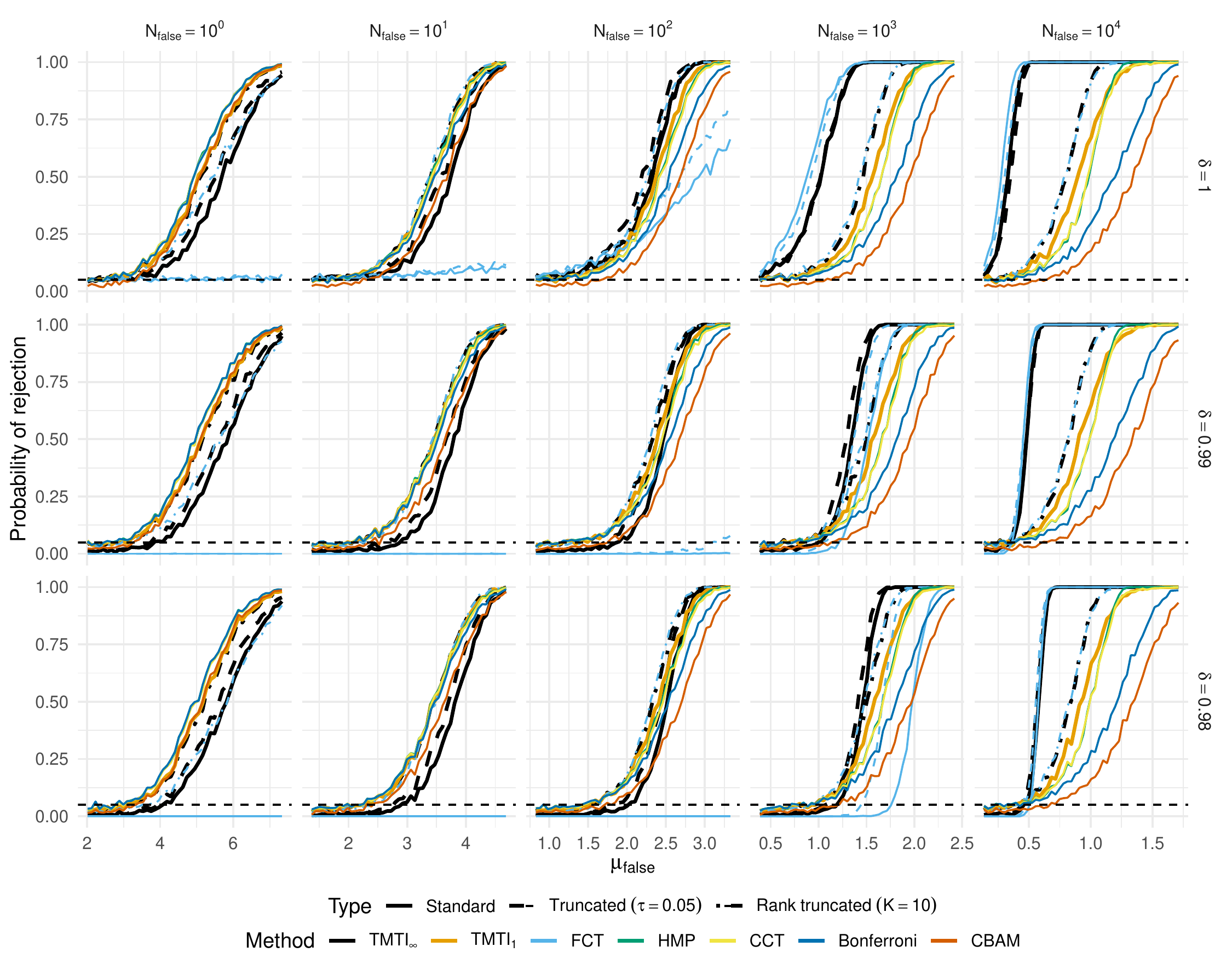}
    \caption{Power curves for different TMTI tests and competing methods. Generally, the TMTI$_\infty$ and tTMTI$_\infty$ work well in all settings. The values of $\mu_{\operatorname{false}}$ are chosen equidistantly between the two values, which satisfy that a Bonferroni test has either 5\% or 99\% power to reject the global null hypothesis in a setting with no conservatism.}
    \label{fig:new_simulation}
\end{figure}
The results of the simulations are displayed in \cref{fig:new_simulation}. Overall, there are three things to notice.

First, the TMTI$_\infty$ and tTMTI$_\infty$ generally
work well no matter how many false hypotheses there are. When there is only a single false hypothesis, these methods have less power than, e.g., a Bonferroni correction, which has the highest power in this scenario, but both methods have considerably higher power than both the {\FCT} and {\TPM}. When there are more false hypotheses, TMTI$_\infty$ and tTMTI$_\infty$ perform on par with the {\FCT} and {\TPM}, having considerably higher power than the remaining methods.
No other methods exhibit this property; the {\CCT}, {\HMP}, {\CBAM} and Bonferroni test all work well when signals are sparse, but have low power when there are many weak signals. In contrast, the {\FCT} and {\TPM}  work well when signals are dense and weak, but have almost no power when signals are sparse and strong. Thus, the TMTI$_\infty$ and tTMTI$_\infty$ appear to have high power against all alternative hypotheses, mimicking the properties of, e.g., the ACAT-O, a test which is shown to have high power against both sparse and dense alternatives \citep{liu2019acat}.
The ACAT-O, however, is designed specifically for sequencing studies and works by leveraging information about the minor-allele counts of a sequencing study, and thus cannot be directly applied in other settings. In contrast, TMTI$_\infty$ and tTMTI$_\infty$ work as regular combination tests and can be applied to any type of data, given that the assumption of independence is satisfied. 

Second, TMTI$_1$ and rtTMTI$_\infty$ work well when signals are sparse, although not better than a Bonferroni correction. When signals are dense, these methods have less power than the TMTI$_\infty$, tTMTI$_\infty$, {\FCT}   and {\TPM}, but higher power than the {\CCT}  and the {\HMP}. The TMTI$_1$ and {\rTPM}  have almost identical power in all settings.

Third, all methods are, in some degree, affected by conservatism, in the sense that all methods generally have less power when the $p$-values from true hypotheses are conservative. When there are few false hypotheses ($N_{\operatorname{false}} \leq 10^2$), the {\FCT} and {\TPM}  have almost no power even under mild conservatism. When there are sufficiently many false hypotheses ($N_{\operatorname{false}} = 10^4$), the effect of conservatism is less pronounced. Overall, it appears that the TMTI tests are less affected by conservatism than the Fisher tests.

It is in line with the intuition behind the TMTI tests that TMTI$_\infty$ and tTMTI$_\infty$ do not perform as well as its competitors in situations where signals are sparse but strong, because these achieve their power from ‘too many’ of the marginal hypotheses being false. In contrast, minimum-$p$ based tests, such as the Bonferroni procedure, need only a single, very strong signal to detect that the global null is false. Similarly, \citet{acat_original} argue that the {\CCT}  only makes use of the first few small $p$-values to represent the overall significance. The same holds true for the {\HMP} \citep{wilson2019harmonic}. The reason that TMTI$_1$ and rtTMTI$_\infty$ still perform similarly to these three methods is that the first local minimum of the sorted and transformed $p$-values is likely to lie early on when there are only a few, very small $p$-values, and likely to coincide with the global minimum of the $K$ smallest $p$-values, if $K$ is sufficiently small. Thus, TMTI$_1$ and rtTMTI$_\infty$ share the property, that they are influenced the most by a few of the smallest $p$-values. The global minimum, however, need not lie early on, when there are only a few false hypotheses, meaning that the few signals that we do observe can potentially be missed when assessing the overall significance using TMTI$_\infty$ or tTMTI$_\infty$.

In \ref{appendix:mixed_betas} we repeated the experiment of this section in a setting with non-constant $\mu$ values (i.e., when $N_{\operatorname{false}} > 1$, the values of $\mu_{\operatorname{false}}$ were allowed to differ for each false marginal hypothesis), finding results similar to those shown in \cref{fig:new_simulation}.
}

\section{Multiple testing and strong FWER control}\label{sec:MTP}
In this section, we consider a common task in statistics. Given $p$-values for a collection of hypotheses, which hypotheses can safely be rejected? As each $p$-value gives marginal Type I error control by definition, a naive approach would be to set a level, $\alpha$, and reject any hypothesis if its corresponding $p$-value falls below $\alpha$. However, as the number of tests conducted increases, more Type I errors will be made, which makes it necessary to employ methods that control for multiple testing. Popular targets one may wish to control for include the False Discovery Rate \citep{benjamini1995controlling} and the Family-Wise Error Rate (FWER). To control the FWER the Bonferroni correction is often used, as it is easy to implement and guarantees strong FWER control. However, the Bonferroni correction has received criticism for, among other things, heavily increasing the risk of making Type II errors, i.e., failing to reject false hypotheses \citep{perneger1998s}. A general approach for turning global testing procedures into a procedure that controls the FWER for elementary hypotheses is the  {\CTP} of \citet{marcus1976closed}. We briefly review the theory on  {\CTPs}.
\begin{definition}\label{def:CTP}
	Let $\mathcal{I}$ be a set of indices and let $\{H_i\}_{i \in \mathcal{I}}$ denote a collection of hypotheses.
	For any subset $\mathcal{J} \subseteq I$, let $H_\mathcal{J} \coloneqq \bigcap_{j \in \mathcal{J}} H_j$ be the joint hypothesis. Let	
	$\phi^\mathcal{J}$ be a random variable on $[0, 1]$ satisfying
	\begin{displaymath}
	H_\mathcal{J} \text{ true } \implies  		\forall \alpha \in [0, 1]: \hspace{1cm} \prop(\phi^\mathcal{J} \leq \alpha) \leq \alpha.
	\end{displaymath}
	That is, $\phi^\mathcal{J}$ is a valid $p$-value for the test of the joint hypothesis $H_\mathcal{J}$.
	Furthermore, for any subset $\mathcal{J}\subseteq \mathcal{I}$ we define the \textbf{closure} of $\mathcal{J}$ in $\mathcal{I}$ to be
	\begin{displaymath}
  	\mathcal{J}_\mathcal{I}^* \coloneqq \bigcup_{{\mathcal{K}}: {\mathcal{K}} \subseteq \mathcal{I},~ \mathcal{J} \subseteq {\mathcal{K}}} \{{\mathcal{K}}\}.
	\end{displaymath}
	That is, $\mathcal{J}_\mathcal{I}^*$ is the set of all supersets of $\mathcal{J}$ that are contained in $\mathcal{I}$.
	A \textbf{Closed Testing Procedure} for the test of the joint hypothesis $H_\mathcal{J}$ is one that rejects $H_\mathcal{J}$ at level $\alpha$ if and only if every superset of $\mathcal{J}$ in $\mathcal{I}$ is also rejected at level $\alpha$, i.e.,
	\begin{displaymath}
  		(H_\mathcal{J} \text{ is rejected at level } \alpha) = \bigcap_{{\mathcal{K}} \in \mathcal{J}_\mathcal{I}^*} \left(
  		\phi^{\mathcal{K}} \leq \alpha
  		\right).
	\end{displaymath}
	That is, the event that we reject $H_\mathcal{J}$ occurs if and only if we reject all supersets of $\mathcal{J}$ in $\mathcal{I}$ marginally.
\end{definition}
From the above, we see that a  {\CTP} is more strict than marginal testing. That is, it becomes more difficult to reject any hypothesis, as we now need to reject all hypotheses that include the hypothesis of interest -- not only the hypothesis itself. The upside is that we obtain strong control of the FWER.
\begin{theorem}[\cite{marcus1976closed}]\label{thm:FWER}
	Let $\mathcal{J}_1, \dots, \mathcal{J}_m$ be distinct subsets of a larger set of indices $\mathcal{I}$. Testing $H_{\mathcal{J}_1}$ through $H_{\mathcal{J}_m}$ each at level $\alpha$ by means of a closed testing procedure controls the FWER at level $\alpha$ in the strong sense.
\end{theorem}
Given any general method to construct tests of joint hypotheses from marginal tests, we can employ these in a  {\CTP} to obtain strict control of the FWER. It is generally accepted that  {\CTP}s~are more powerful than other methods that control the FWER \citep{grechanovsky1999closed}, although this power increase comes at the cost of a heavy computational burden. Given $m$ marginal hypotheses which we want to test, we need to perform $\sum_{i=1}^m \binom{m}{i} = 2^m - 1$ tests. This is because we need to test all possible intersection hypotheses, which corresponds to the powerset of all hypotheses, minus the empty set. Thus, in many cases, it is not feasible to perform a  {\CTP} when $m$ is even slightly large. Indeed, even with just $m=300$ marginal tests, the number of tests to be performed in a full  {\CTP} is $2^{300} - 1 \approx 2 \cdot 10^{90}$ -- roughly $10$ billion times the number of atoms in the observable universe. Thus, with many procedures, one seeks to find a shortcut so that only a subset of the powerset of hypotheses needs to be tested. This is often possible \citep{grechanovsky1999closed} and considerably reduces the computational complexity of carrying out a  {\CTP}.

\citet{zaykin2002truncated} introduce a shortcut for the {\TPM}, reducing the computational complexity of the  {\CTP} from $\mathcal O(2^m)$ to $\mathcal O(m^2)$.
In a recent result, \citet{tian2021largescale} give the same shortcut for a family of combination tests that are sums of marginal tests. \changed{%
\citet{dobriban2020fast} gives a shortcut for test statistics that are monotone and symmetric. Here, we provide a shortcut for class of combination tests, which are monotone but not necessarily symmetric, and not necessarily sums of marginal tests. Furthermore, we show that TMTI$_\infty$, tTMTI$_\infty$ and rtTMTI$_\infty$ all admit this shortcut. 
}

\begin{lemma}\label{lemma:minSubset}
	Let $\boldsymbol{p}_\mathcal{I}$ be a set of observed $p$-values with $\mathcal I \coloneqq \{1, \dots, m\}$ and $m \geq 2$. Let $\mathcal J^k$ be the set of all subsets of $\mathcal{I}$ with $\vert\mathcal{J}\vert = k$. 
	Let $\mathcal X \subseteq \mathbb{R}$ be a set and let $F_{(1)}: [0,1] \to \mathcal{X}, \dots, F_{(k)}: [0, 1] \to \mathcal{X}$ be a sequence of functions that satisfy
	\begin{align*}
    \forall j \in \mathcal{I}~~ \forall x \in \mathcal{X} ~~ \forall \epsilon \geq 0: 
    & \hspace{1cm}
	F_{(j)}(x) \leq F_{(j)}(x + \epsilon)\label{eq:assumption_weak}\tag{C1}
	\shortintertext{and}
	\forall j \in \{1, \dots, m-1\} ~~ \forall x \in \mathcal X : 
	& \hspace{1cm}
	F_{(j)}(x) \geq F_{(j + 1)}(x).\label{eq:assumption_F}\tag{C2}
	\end{align*}
Define for all $\mathcal J \in \mathcal{J}^k$ the random variable $\boldsymbol{Y}^\mathcal{J} \coloneqq (F_{(1)}(p_{(1)}), \dots, F_{(k)}(p_{(k)}))$
	and let $h: \mathcal{X}^k \to [0,1]$ be a function satisfying
	\begin{equation}\label{eq:assumption_h}\tag{C3}
\forall \boldsymbol{x} \in \mathcal{X}^k ~~ \forall \boldsymbol{\epsilon} \in \mathbb{R}_{+}^k: \hspace{1cm}
		h(\boldsymbol{x}) \leq h(\boldsymbol{x} + \boldsymbol{\epsilon}),
	\end{equation}
Let $\eta:\mathcal I \to \mathcal I$ be a bijection ordering $\boldsymbol{p}_\mathcal{I}$, i.e., $p_{\eta(1)} \leq \dots \leq p_{\eta(m)}$.
	It then follows that
	for any two sets, $\mathcal J_1, \mathcal J_2 \in \mathcal J^k$
	\[%
	\eta(\mathcal{J}_1) \leq \eta(\mathcal{J}_2) \implies h(\boldsymbol{Y}^{ \mathcal J_1}) \leq h(\boldsymbol{Y}^{ \mathcal J_2}).
	\]
	In the above, the operation $\leq$ applied to the sets $\eta(\mathcal{J}_1)$ and $\eta(\mathcal{J}_2)$ is taken to mean element-wise less than or equal to.
\end{lemma}
Lemma \ref{lemma:minSubset} states that whenever we consider $k$ $p$-values, we will obtain a smaller test statistic if we substitute one or more of them with smaller $p$-values. In the context of closed testing, this implies that when considering the closure of an atom, say $\{j\}$, then among all subsets of size $k$ in $\{j\}_{\mathcal{I}}^*$ we do not need to test all $m! / ((m-k)!k!)$ intersection hypotheses. This is because we know that the largest (smallest) test statistic is obtained by considering the $p$-value $p_j$ combined with the $k-1$ largest (smallest) remaining $p$-values. Assuming that the underlying distribution of the $p$-values is exchangeable, this implies that the $p$-values for the combination tests obey the same inequalities as the test statistics, and thus we need only consider the intersection hypothesis which we know will yield the largest $p$-value. 

\begin{remark}\label{remark:minSubset}
The same result as in Lemma \ref{lemma:minSubset} can be obtained by reversing the inequalities in Equations \eqref{eq:assumption_weak}, \eqref{eq:assumption_F} and \eqref{eq:assumption_h}. In contrast, we can obtain a version which gives
    \[
    \eta(\mathcal{J}_1) \leq \eta(\mathcal{J}_2) 
    \implies
    h(\boldsymbol{Y}^{\mathcal{J}_1}) \geq h(\boldsymbol{Y}^{\mathcal{J}_2})
    \]
if we reverse the inequalities in Equations \eqref{eq:assumption_weak} and \eqref{eq:assumption_F} and keep Equation \eqref{eq:assumption_h}, or if we reverse the inequality in Equation \eqref{eq:assumption_h} and keep Equations \eqref{eq:assumption_weak} and \ref{eq:assumption_F}. The choice of which version to use depends on whether small or large values of the test statistic are critical.
\end{remark}

\begin{theorem}\label{theorem:main0}
	TMTI$_\infty$, tTMTI$_\infty$ and rtTMTI$_\infty$ all satisfy the conditions of Lemma \ref{lemma:minSubset}.
\end{theorem}
\begin{remark}
	Even though the TMTI$_\infty$ variants all satisfy the conditions in Lemma \ref{lemma:minSubset}, not all TMTI$_n$ variants do. To see this, consider two sets of $p$-values, $\boldsymbol{p}_1 = (0.25, 0.50, 0.75)$ and $\boldsymbol{p}_2 = (0.2, 0.5, 0.75)$. Then $\boldsymbol{Y}_1 = (0.58, 0.5, 0.42)$, making $0.42$ the first local minimum of $\boldsymbol{Y}_1$, and $\boldsymbol{Y}_2 = (0.49, 0.5, 0.42)$, making $0.49$ the first local minimum  of $\boldsymbol{Y}_2$. Thus, we have $\boldsymbol{Y}_2 \leq \boldsymbol{Y}_1$ but $h_{TMTI_1}(\boldsymbol{Y_2}) > h_{TMTI_1}(\boldsymbol{Y_1})$.
\end{remark}

\begin{theorem}\label{theorem:main}
	Let $p_\mathcal{I}$, $F_{(1)}, \dots, F_{\vert\mathcal I\vert}$, $\mathcal X$, $h$ and $\boldsymbol{Y}^\mathcal{J}$ be defined as in Lemma \ref{lemma:minSubset}. Assume that the underlying distribution of $p_\mathcal{I}$ is exchangeable. If Equations \eqref{eq:assumption_weak}, \eqref{eq:assumption_F} and \eqref{eq:assumption_h} are satisfied, then a  {\CTP} using $h(\boldsymbol{Y}^\mathcal{J})$ as test statistic can be used to obtain control of the FWER for all marginal hypotheses in at most $\frac{1}{2} m (m - 1)$ steps.
\end{theorem}

\begin{remark}
	The result in Lemma \ref{lemma:minSubset} and its converse in Remark \ref{remark:minSubset} does not only apply to TMTI statistics. For example, letting $F_{(1)}, \dots, F_{(m)}$ be the identity mappings and $h(\boldsymbol{x}) \coloneqq - 2\sum_{i = 1}^m \log x_i$ we obtain the {\FCT},  which then gives us the well-known $\mathcal{O}(m^2)$ shortcut described e.g.,in \citet{zaykin2002truncated}. Similarly, letting $h(\boldsymbol{x}) = \frac{1}{m} \sum_{i=1}^m \tan((0.5-x_i)\pi)$, we obtain the unweighted {\CCT}.
\end{remark}

In Algorithm \ref{alg:CTP}, we give an example of how this shortcut procedure can be implemented to return adjusted $p$-values for the tests of all marginal hypotheses.

\begin{algorithm}
\caption{Shortcut  {\CTP} for statistics satisfying Conditions \eqref{eq:assumption_weak}, \eqref{eq:assumption_F} and \eqref{eq:assumption_h}}
\label{alg:CTP}
    \DontPrintSemicolon
    \KwIn{Sorted $p$-values $p_1 \leq \dots \leq p_m$ for tests of hypotheses $H_1, \dots, H_m$, a significance level $\alpha$}
	\KwOut{Adjusted $p$-values for the tests of $H_1, \dots, H_m$}
	
	Construct an empty $m \times m$ matrix $Q$.
	
	\For{%
	    $i = 1, \dots, m - 1$
	}{%
	    $c \leftarrow m$.
	    
	    $Q_{i, c} \leftarrow p_i$.
	    
    	\For{%
            $j = m, \dots, i + 1$
        }{%
            $c \leftarrow c - 1$.
        
            Test the hypothesis $\left(\bigcap_{k = m}^{j} H_k\right) \cap H_i$ and save the $p$-value as $p_{i, m:j}$.
            
            $Q_{i, c} \leftarrow p_{i, m:j}$.
        }
	}
	
	\For{%
	    $i = 1, \dots, m - 1$
	}{%
	    $Q_{i, (i + 1):m} \leftarrow Q_{i, i}$
	    
	    $\tilde{p}_i \leftarrow \max Q_{1:m, i}$
	}
	$\tilde{p}_m \leftarrow  \max Q_{1:m, m}$
	
	Return $\tilde{p_1}, \dots, \tilde{p}_m$ as adjusted $p$-values for $H_1, \dots, H_m$.
\end{algorithm}

If the practitioner is content with having only a lower bound on the adjusted $p$-value, whenever the test is not rejected at a chosen level, Theorem \ref{theorem:main} provides an upper bound on the number of steps required to complete the  {\CTP}. For instance, if the global hypothesis cannot be rejected, no marginal hypothesis can be rejected, and the procedure can therefore be stopped and the $p$-value for the test of the global hypothesis can be used as a lower bound for all adjusted $p$-values. Stopping the procedure early can speed up computations considerably, especially when $m$ is large but very few hypotheses can be rejected. \changed{%
Furthermore, \cref{lemma:minSubset} also implies that an adjusted $p$-value for a single elementary hypothesis can be obtained in only $m$ steps. In practice, one is often only interested in obtaining adjusted $p$-values for the hypotheses, for which the marginal $p$-value is significant (as the remaining hypotheses cannot be rejected in a {\CTP}). Thus, if there are, say, $n$ elementary hypotheses for which the marginal $p$-value is significant, these can be adjusted in $\mathcal{O}(nm)$ complexity. 
}

\begin{figure}[hbt]
	\centering
\begin{tikzcd}
    &                      &                        & {H_{1,2,3,4}} \arrow[ld] \arrow[rd] &                                   &                                 &     \\
    &                      & {H_{1,3,4}} \arrow[ld] &                                     & {H_{2,3,4}} \arrow[ld] \arrow[rd] &                                 &     \\
    & {H_{1,4}} \arrow[ld] &                        & {H_{2,4}} \arrow[ld]                &                                   & {H_{3,4}} \arrow[ld] \arrow[rd] &     \\
H_1 &                      & H_2                    &                                     & H_3                               &                                 & H_4
\end{tikzcd}
  \caption{Test procedure with $m=4$}
  \label{fig:testtree}
\end{figure}
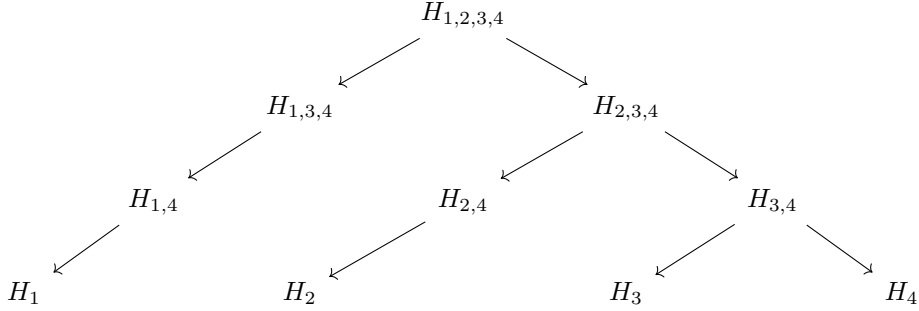

We have exemplified in \cref{fig:testtree} how the reduced test-tree looks for the  {\CTP} when applying any test that satisfies the conditions of Lemma \ref{lemma:minSubset} in the case of $m = 4$ total tests, where the hypotheses are labeled such that $H_i$ corresponds to the $i$th lowest $p$-value. To obtain an adjusted $p$-value for any marginal hypothesis, say $H_i$, one takes the maximal $p$-value from the test of all ancestral hypotheses in the graph. For example, the adjusted $p$-value for the test of $H_2$ would be the largest of the $p$-values from the tests of $H_2$, $H_{2,4}$, $H_{2,3,4}$ and $H_{1,2,3,4}$.

\subsection{A remark on mixture strategies in Closed Testing Procedures}\label{sec:mix}
Definition \ref{def:CTP} of a  {\CTP} and the subsequent Theorem \ref{thm:FWER} on FWER control make no assumptions on the choice of local tests, $\phi^\mathcal{J}$, and these can in principle vary across all intersection hypotheses to be tested. We only require that every local test, $\phi^\mathcal{J}$, is a valid $\alpha$-level test. The natural choice is to use the same kind of test at every intersection hypothesis, e.g.,  TMTI$_\infty$. However, we can in principle employ any choice of local test. In some cases, we argue, it is reasonable to use different local tests for different kinds of intersection hypotheses. When we go through a  {\CTP}, we are going to consider tests of many different kinds of hypotheses, and in particular different kinds of alternative hypotheses. As previously discussed, TMTI$_\infty$ has slightly lower power than other methods in situations where signals are extremely sparse. According to the shortcut strategy outlined in Algorithm \ref{alg:CTP}, we need only consider the $\vert\mathcal J\vert - 1$ largest $p$-values alongside the $j$th $p$-value, when testing all supersets of $\{j\}$ of size $\vert\mathcal J\vert$. When $\vert\mathcal J\vert$ is small or when false hypotheses are sparse, it is likely that these intersection hypotheses will consist of a single false hypothesis (if any) with all the remaining hypotheses being true. It makes sense, then, to employ a different test in these situations, which has power against alternatives of sparse signals. However, this alternative test needs to satisfy the same shortcut as TMTI$_\infty$ in order to be employed across all supersets of equal size. One such choice is rtTMTI$_\infty$ with a low choice of $K$ -- e.g., $K=1$ -- as this method satisfies the same shortcut as TMTI$_\infty$, but has higher power against sparse alternatives, as discussed in \cref{sec:newpower}. Given a number of hypotheses, say $m$, of which we expect $F < m$ to be false, we could for example conduct the  {\CTP} by employing rtTMTI$_\infty$ with a small $K$ whenever we consider supersets of size at most $m - F$ and TMTI$_\infty$ whenever we consider supersets of size greater than $m - F$. We call such a strategy a mixture  {\CTP}. We return to mixture  {\CTPs} in \cref{sec:realdata}, where we analyse a real dataset. The reasoning is that once we start considering supersets of size greater than $m - F$, the intersections hypotheses considered in the shortcut procedure will potentially include multiple false hypotheses, while they will include at most one false hypothesis when considering supersets of size less than $m - F$.

We stress that when employing a mixture  {\CTP}, the choice of local tests should be made \textit{a priori} and not be data-driven, so as not to incur new multiplicity problems.

\section{The number of false hypotheses in a rejection set and $k$-FWER control}\label{sec:quantify}
Given a set of hypotheses, indexed by $\mathcal J$, such that we can safely reject the joint hypothesis $H_\mathcal{J}$ -- i.e., we conclude that at least one of the hypotheses in $\mathcal J$ is false -- the natural question is then how many of the hypotheses in $\mathcal J$ are false. To answer this, \citet{goeman2011multiple} provide a simple way of generating a $1-\alpha$ confidence set for the number of false hypotheses contained in $\mathcal J$ when using  {\CTP}s. Let $\mathcal R \subseteq \operatorname{powerset}(\mathcal J)\backslash \{\emptyset\}$ be the set of all intersection hypotheses that can be rejected by any  {\CTP}. Define $\tau$ to be the number of true hypotheses in $\mathcal J$ and
		$
		t_\alpha \coloneqq \max\{\vert\mathcal K\vert: \mathcal K \subseteq \mathcal J, \mathcal K \not\in \mathcal R \}
		$
	to be the size of the largest intersection hypothesis in $\mathcal J$ that can not be rejected by the  {\CTP}.
	\begin{theorem}[\citet{goeman2011multiple}]\label{thm:goeman}
		The sets
			$
			\{0, \dots, t_\alpha\}
			$
		and
			$
			\{\vert\mathcal J\vert - t_\alpha, \dots, \vert\mathcal J\vert\}
			$
		are $1-\alpha$ confidence sets for the number of true hypotheses, $\tau$, and the number of false hypotheses, $\vert\mathcal J\vert-\tau$, respectively. 
	\end{theorem}
This remarkably simple theorem has the implication, that we can generate confidence sets for the number of false hypotheses among all those tested in only $t_\alpha$ steps when using any test procedure satisfying the conditions of Lemma \ref{lemma:minSubset}, assuming that the $p$-values are realized from an exchangeable distribution. We describe in Algorithm \ref{alg:confSet3} how to do this. 

The quantity $t_\alpha$ depends on the choice of test used on each intersection. Different tests have power against different alternatives, and a test that has low power for a particular intersection hypothesis will be more likely to not reject that hypothesis. Thus, if the chosen test has low power for the particular data, the resulting confidence set for the number of false hypotheses will be conservative. In contrast, if the chosen test has high power against the particular alternative, the confidence set tightens.

As noted in \citet{goeman2011multiple}, we can also apply Theorem \ref{thm:goeman} as a way of controlling the $k$-FWER, i.e., the probability of making at least $k$ false rejections. To control the $k$-FWER, we find the largest $i$ such that $t_\alpha < k$ when calculated for the set of hypotheses yielding the $i$ smallest $p$-values. That is, if we find that $F$ hypotheses in a set, say $\mathcal J$, are false with $1-\alpha$ confidence, then we can reject every hypothesis in that set while controlling the $k$-FWER at $k = \vert\mathcal J\vert - F + 1$. This can be done in $\mathcal O(m^3)$ time. This is because determining $t_\alpha$ from a set $\mathcal J$ is done in $\mathcal O(\vert\mathcal J\vert^2)$ time, as described in Algorithm \ref{alg:confSet3}, and we now need to do this for subsets of increasing (or decreasing, depending on the search direction) size.
It is worth noting, that the $k$ at which the practitioner wishes to control the $k$-FWER need not be chosen \textit{a priori}, as the confidence sets are simultaneous for all choices of $\mathcal J$, simply because the closure of each $\mathcal J$ is contained within the full closure.

\begin{algorithm}
	\caption{Confidence set for the number of false hypotheses among $\mathcal{J} \subseteq \mathcal{I}$}
	\label{alg:confSet3}
	\KwIn{%
		Hypotheses $(H_i)_{i \in \mathcal{I}}$, an ordered set $\mathcal J \subseteq \mathcal I$, ordered $p$-values $p_{1} \leq \dots \leq p_{\vert \mathcal{I} \vert}$ .
	}
	\KwOut{%
	A $1-\alpha$ confidence set for the number of false hypotheses in $\mathcal J$.
	}
	
	\init{}{%
	
	$t_\alpha \leftarrow 0$
	
	$m \leftarrow \vert\mathcal J\vert$
	
	$\tilde{p} \leftarrow (p_j)_{j \in \mathcal J}$
	}
	
	\If{%
		$m = \vert \mathcal{I} \vert$
	}{%
		\For{%
			$i = m, \dots, 1$
		}{%
		$\tilde{\mathcal{J}} \leftarrow \{m- i + 1, \dots, m\} $

		Let $p_{\tilde{\mathcal{J}}}$ be the $p$-value for the test of $\bigcap_{j \in \tilde{\mathcal{J}}} H_{j}$.
	
			\If{%
				$p_{\tilde{\mathcal{J}}} \geq \alpha$
			}{%
				Set $t_\alpha \leftarrow i$
		
				Break the loop and return $\{\vert\mathcal{J}\vert - t_\alpha, \dots, \vert\mathcal{J}\vert\}$.
			}
		}
	}\Else{%
		\For{%
			$i = m ,\dots, 1$
		}{%
		$\tilde{\mathcal J} \leftarrow \{m- i + 1, \dots, m\} $

		$\hat{p} \leftarrow (p_i)_{i \in  \mathcal I \backslash \tilde{\mathcal J}}$ %
    
	    	Let $p_{\tilde{\mathcal J}}$ be the $p$-value for the test of $ \bigcap_{j \in {\tilde{\mathcal J}}} H_{j}$
    
		\For{%
			$j = 1, \dots, \vert \mathcal I \backslash \tilde{\mathcal J} \vert$
		}{%
			Define ${\tilde{\mathcal J}}_2$ as $\tilde{\mathcal J}$ appended with the $j$ largest values of $\hat{\mathcal{I}}$.

			Update $p_{\tilde{\mathcal J}}$ as the $p$-value from the test of $\bigcap_{j \in {\tilde{\mathcal J}}_2} H_{j}$.
		
			\If{%
				$p_{\tilde{\mathcal J}} \geq \alpha$
			}{%
				Set $t_\alpha \leftarrow i$
		
				Break the loop and return $\{\vert\mathcal J\vert - t_\alpha, \dots, \vert\mathcal J\vert\}$.

			}
		}
	}
	}
	
	Return $\{\vert\mathcal{J}\vert - t_\alpha, \dots, \vert\mathcal{J}\vert\}$.
\end{algorithm}

\section{Additional computational considerations}
The computational time of a carrying out a  {\CTP} when using the above shortcuts is manageable for reasonable values of $m$. For example, computing adjusted $p$-values for a set of $m = 100$ $p$-values take roughly two seconds on a standard laptop using single-threaded computations.
Still, there is a considerable amount of computational effort involved in carrying out a  {\CTP} when $m$ is sufficiently large, in part because the CDFs of the TMTI statistics will have to be bootstrapped. To further reduce the computational burden, we offer the following result.
\begin{lemma}\label{lemma:skipM}
    Let $\mathcal{J}_1$ and $\mathcal{J}_2$ be sets such that $\mathcal{J}_1 \subsetneq \mathcal{J}_2$.
    Then
    \[
    \forall x \in (0,1): \hspace{1cm}
    \gamma^{\mathcal{J}_1}(x) < \gamma^{\mathcal{J}_2}(x).
    \]
    That is, a conservative $p$-value for the test of $H_{\mathcal{J}_1}$ can be obtained by using $\gamma^{\mathcal{J}_2}$ instead of $\gamma^{\mathcal{J}_1}$ when computing the $p$-value.
\end{lemma}
The purpose of Lemma \ref{lemma:skipM} is that the user can choose to skip the bootstrap at several layers of the  {\CTP} and instead simply use the CDF of a higher layer. This improves the running time of the algorithm at the cost of using conservative $p$-values at the layers where the bootstrap was skipped. Exactly how costly this trade-off is, depends on how many layers are skipped each time. We conjecture that the $p$-value will only be slightly conservative if the number of layers skipped is small relative to the size of the subsets considered.

\section{An example using real data}\label{sec:realdata}
In this section, we give an example of how TMTI$_\infty$ performs against other methods when applied to real data. For this purpose, we consider data from the \textit{National Assessment of Educational Progress} on the state-wise changes in eighth-grade mathematics achievements from 1990 to 1992. This data is, among other places, presented in \citet{williams1999controlling}, where the authors compute two-sided $T$-tests for the mean change in mathematical achievements over the two-year period, to quantify whether or not any particular state has progressed or worsened during that time period. The original data includes one $p$-value of exactly $0$. We have rounded this to be $0.00001$ here to ensure Type I error control. The same data is used in \citet{benjamini2000adaptive}, where the authors find that mathematics achievements have changed significantly in $24$ of the $34$ states. However, the authors control the more lenient False Discovery Rate, and thus their results are not directly comparable to the ones presented here. \citet{williams1999controlling} apply a significance level of $0.025$ instead of the usual $0.05$. In \citet{benjamini2000adaptive}, the authors apply the usual $0.05$ significance level but have doubled all $p$-values such that the results are comparable. We have done the same here. The data, adapted from \citet{williams1999controlling}, is presented in Table \ref{table:realdata}. Given that we only have access to summary statistics, we assume that all $p$-values are independent. Whether this is a reasonable assumption can be debated.

There are several questions regarding this data that may be of interest to the practitioner:
\begin{enumerate}
	\item Did mathematics achievements change significantly in any state from 1990 to 1992?
	\item In how many states did mathematics achievements change significantly?
 	\item In what states did mathematics achievements change significantly?
\end{enumerate}
To answer the first question, we can, for example, apply TMTI$_\infty$ to obtain a $p$-value. Doing so results in a $p$-value of $1.58\times 10^{-13}$. Thus, we find evidence that mathematics education has changed significantly in at least one of the $34$ states.

To answer the remaining two questions, we 
apply TMTI$_\infty$ in a  {\CTP} as well as a mixture  {\CTP}, using rtTMTI$_\infty$ with $K = 1$ whenever we consider fewer than $15$ hypotheses and TMTI$_\infty$ when considering more. For comparison, we also apply the standard Bonferroni correction as well as the {\FCT}  in a  {\CTP} and a {\rTPM}/{\FCT} (denoted rTPM/FCT)  mixture, using the {\rTPM} with $K = 1$ for intersection hypotheses of size at most $15$. The choice of $K = 1$ for intersection hypotheses smaller than $15$ corresponds to a belief that at most $15$ hypotheses are true. Here, we have chosen $15$ at random, but a practitioner with subject matter knowledge can choose this based on prior knowledge.

By applying Algorithm \ref{alg:confSet3} with TMTI$_\infty$ we find that $\{23, \dots, 34\}$ is a $95\%$ confidence set for the number of states in which mathematics achievements have changed. In contrast, using the {\FCT},  the confidence set is $\{19, \dots, 34\}$. The same confidence set is found when applying the mixture strategies. That is, we can say with $95\%$ confidence, that mathematics achievements have changed significantly in at least $23$ states when using TMTI$_\infty$. The improved performance of TMTI$_\infty$ over the {\FCT} here is likely due to TMTI$_\infty$ generally having high power across a wide range of settings (see \cref{sec:newpower}), whereas the {\FCT} lacks power in settings with sparse, strong signals. When carrying out a {\CTP}, we test many different kinds of joint hypotheses (i.e., some containing more false hypotheses than others), and it is thus beneficial that the employed test has high power in many different settings.

To determine which of the hypotheses we can say with certainty are false while controlling the FWER, we applied Algorithm \ref{alg:CTP}. Here, TMTI$_\infty$, the {\FCT}  and the Bonferroni correction perform identically and are capable of rejecting the bottom four hypotheses. In contrast, when we apply the two mixture strategies, we can reject the bottom seven hypotheses. Thus, by incorporating prior knowledge, we can increase the size of the rejection set for which we have FWER control.

\begin{figure}[H]
  \centering
  \includegraphics[width=7.5cm, height=7.5cm]{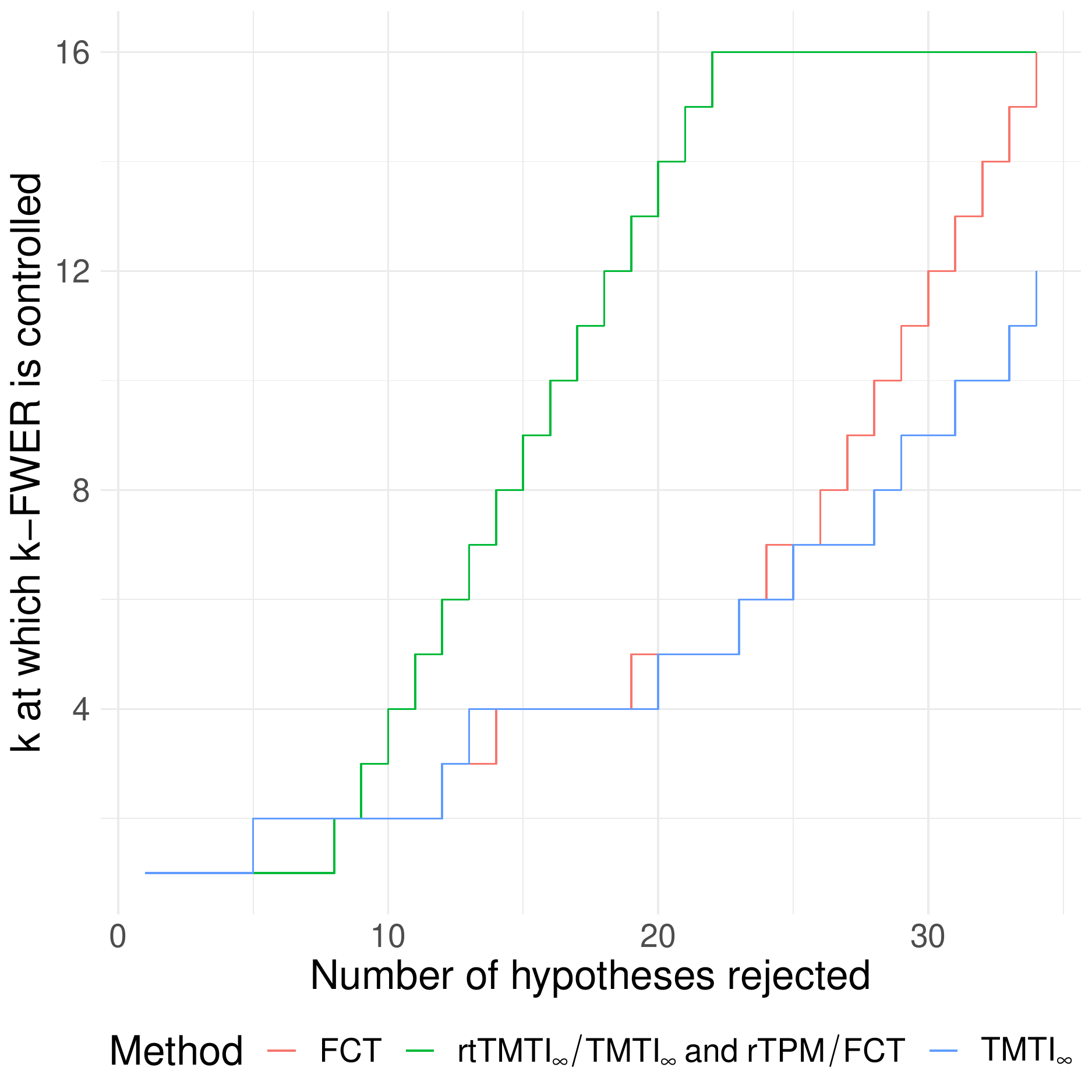}
  \caption{Overview of the different  {\CTP} methods employed and the $k$ at which they control the $k$-FWER, when rejecting the bottom $t$ hypotheses, for $t = 1,\dots, 34$. The two mixture strategies, rtTMTI$_\infty$/TMTI$_\infty$ and rTPM/FCT,  are colored the same, as their results are identical.}
  \label{fig:kFWERresults}
\end{figure}

To find which of the hypotheses can be rejected while controlling the more general $k$-FWER, we consider sets of increasing size of the smallest $p$-values, each time calculating $t_\alpha$. For any chosen set, say $\mathcal{J}$, we can reject the entire set while controlling the $k$-FWER at $t_{\alpha} + 1$. Doing so, we find that the 11 hypotheses giving rise to the smallest $p$-values can be rejected by the TMTI$_\infty$ while controlling the $k$-FWER at $k=2$. In contrast, the mixture  {\CTP} is only capable of rejecting the bottom eight hypotheses while controlling the $k$-FWER at $k=2$. If we are willing to accept a more lenient $k$-FWER control of $k=5$, we are capable of rejecting the bottom 22 hypotheses using the TMTI$_\infty$ and the bottom 11 hypotheses using the rtTMTI$_\infty$/TMTI$_\infty$ mixture. In \cref{fig:kFWERresults}, we have summarized the associated $k$ at which we control the $k$-FWER at, when rejecting the bottom $t$ hypotheses, for $t=1,\dots, 34$, when using TMTI$_\infty$, the rtTMTI$_\infty$/TMTI$_\infty$ mixture, the {\FCT} and the rTPM/FCT  mixture, respectively. Here, we note that TMTI$_\infty$ is weakly better than the {\FCT}   except for at a single rejection set, $\{1, \dots, 13\}$.

\begin{table}[hbt]
\begin{tabular}{l|llllllll}
State 		& GA & AR & AL & NJ & NE & ND & DE & MI \\
$p$-value (\%)   & 85.628 & 60.282 & 44.008 & 41.998 & 38.640 & 36.890 & 31.162 & 23.522 \\
                 & LA & IN & WI & VA & WV & MD & CA & OH \\
                 & 20.964 & 19.388 & 15.872 & 14.374 & 10.026 & 8.226 & 7.912 & 6.590 \\
                 & NY & PA & FL & WY & NM & CT & OK & KY \\
                 & 5.802 & 5.572 & 5.490 & 4.678 & 4.650 & 4.104 & 2.036 & 0.964 \\
                 & AZ & ID & TX & CO & IA & NH & NC & HI \\
                 & 0.904 & 0.748 & 0.404 & 0.282 & 0.200 & 0.180 & 0.002 & 0.002 \\
                 & MN & RI &    &    &    &    &    & \\
                 & 0.002 & 0.001 &         &         &         &         &         &
\end{tabular}
\caption{States and their $p$-values for $T$-tests of changes in mathematics achievements from 1990 to 1992.}
\label{table:realdata}
\end{table}

That the two mixture strategies have higher power to detect differences when controlling the FWER than TMTI$_\infty$ and the {\FCT}, but lower power to detect differences when using the more lenient $k$-FWER control, may seem counter-intuitive and requires some exposition. The difference lies in what intersection hypotheses need to be considered in the full  {\CTP} test tree. When controlling the FWER, we are in principle looking through all intersection hypotheses, and then we use the maximal $p$-values along the closure of each atom as the adjusted $p$-values. As outlined in Lemma \ref{lemma:minSubset}, however, we need only consider the part of the closure that contains subsets containing only the atom and the largest $p$-values. When we apply Algorithm \ref{alg:confSet3} iteratively to obtain $k$-FWER control, we are considering the closures, not of atoms, but of intersections. Put differently, consider the index set $\mathcal J \coloneqq \{1, \dots, t\}$, for some $t$, of the $t$ smallest $p$-values. We are then going to consider the closure of $\{t\}$ in $\mathcal J$, i.e., $\{t\}_\mathcal{J}^*$ at first. For all sets in this closure, say $\tilde{\mathcal{J}} \in \{t\}_\mathcal{J}^*$, we are then going to calculate the adjusted $p$-value as the maximal $p$-value along the closure of $\tilde{\mathcal{J}}$ in $\mathcal I$, i.e., $\tilde{\mathcal{J}}_\mathcal{I}^*$. These sets are not, as they were in the ordinary  {\CTP}, sets consisting of an atom unioned with the largest $p$-values, but rather several, possibly neighboring, atoms unioned with the largest $p$-values. We constructed the mixture strategies to have higher power in situations in which we considered intersection hypotheses with only a single false hypothesis in them. Using this method, we are now considering sets that possibly have multiple false hypotheses in them, even when the total number of hypotheses included in the set is low -- which is when TMTI$_\infty$ gains its power. In contrast, rtTMTI$_\infty$ with a small $K$ loses power, when there are more than $K$ false hypotheses present in the intersection hypothesis.

A detailed table with adjusted $p$-values for all of the tests employed here can be found in \ref{appendix:pvalue_table}.

\section{Conclusion}
We have introduced the ‘Too Many, Too Improbable’ (TMTI) family of combination test statistics for testing joint hypotheses among $m$ marginal hypotheses. The TMTI family includes truncation-based tests, similar to those of \citet{zaykin2002truncated} and \citet{dudbridge2003rank}, for testing global hypotheses against sparse alternatives. We have shown in \cref{sec:newpower} that the TMTI tests outperforms other combination tests in many situations. 
\changed{In particular, we found that TMTI$_\infty$ and tTMTI$_\infty$ were the only tests that were able to achieve high power both when signals are dense but weak and when signals are sparse but strong. Although we found in all scenarios that there was at least one other test that performed equally as well as TMTI$_\infty$ tTMTI$_\infty$, no other combination tests had similar performance across all scenarios. This property is useful, e.g., if one has no \textit{a priori} knowledge about the sparsity and strength of signals and for generating $1 - \alpha$ confidence sets for the number of false hypotheses in a rejection set.}

In \cref{sec:MTP}, we have given an $\mathcal{O}(m^2)$ shortcut for controlling the Family-Wise Error Rate using  {\CTPs} \citep{marcus1976closed} for a large class of test statistics, which includes the TMTI family of test statistics, but also the {\CCT} among others. Using this shortcut, we use the work of \citet{goeman2011multiple} in \cref{sec:quantify} to develop an $\mathcal{O}(m^3)$ algorithm for controlling the generalized FWER as well as an $\mathcal{O}(m)$ algorithm for obtaining $1-\alpha$ confidence sets for the number of false hypotheses among all hypotheses.

In \cref{sec:realdata}, we applied a TMTI test in a  {\CTP}, as well as a mixture -- i.e., varying local tests across the  {\CTP} -- of two TMTI tests, to a real dataset and compared it to the Fisher Combination Test applied in a  {\CTP}. Here we found that all TMTI tests were able to reject the same hypotheses as the {\FCT},  but that the TMTI test generated a narrower confidence set for the number of false hypotheses among the collection of considered hypotheses. Additionally, we found that by employing mixture strategies, we were able to reject more hypotheses than with standard methods. However, the mixture strategies performed worse when controlling the $k$-FWER with $k\geq 2$.

\section*{Supplementary material}
Proofs of all lemmas and theorems are supplied in the appendix. Furthermore, the appendix includes further simulations to support those of \cref{sec:newpower} and a detailed table containing the adjusted $p$-values of the procedures applied in \cref{sec:realdata}.

The TMTI family of test statistics and the shortcuts for $k$-FWER and confidence sets is implemented in the R package \texttt{TMTI} available on CRAN.

\clearpage
\appendix

\section{Proofs}
\subsection{Proof of Lemma \ref{lemma:dominated}}
	This holds trivially, as the first $Y$ smaller than the following $n$ must necessarily be larger when we only consider the first $c < \vert\mathcal{I}\vert$ values of $\boldsymbol{Y}^\mathcal{I}$ than if we consider the full sequence.

\subsection{Proof of Theorem \ref{theorem:gamma}}
Assume first that $c = m$.
The order statistics $(P_{(1)}, \dots, P_{(m)})$ have a constant joint density $m!$ on the simplex $\{\boldsymbol{p} \in [0,1]^m: p_1 \leq \dots \leq p_m \}$. Thus
    \begin{align*}
    \gamma^{\mathcal{I}}(x) \coloneqq&
    \prop(\min ( \beta(1,m)(P_{(1)}), \dots, \beta(m, 1)(P_{(m)})) \leq x) \\
    =& 1- \prop(\beta(1,m)(P_{(1)}) > x, \dots, \beta(m, 1)(P_{(m)}) > x) \\
    =& 1 - \prop(P_{(1)} > {\beta^{-1}(1, m)(x)}, \dots, P_{(m)} > {\beta^{-1}(m, 1)(x)}) \\
    =& 1 - \prop(P_{(1)} > x_1, \dots, P_{(m)} > x_m) \\
    = & 1 - m! \int_{x_m}^{1} \int_{x_{m-1}}^{q_m} \dotsc \int_{x_{2}}^{q_3} \int_{x_1}^{q_2}~
    dq_1~dq_2 \dotsc  dq_{m-1} ~ dq_m. \numberthis \label{eq:integral}
    \end{align*}
The integral in \cref{eq:integral} can be expressed as %
    \[
    \int_{x_m}^{1} \int_{x_{m-1}}^{q_m} \dots \int_{x_{2}}^{q_3} \int_{x_1}^{q_2}~
    dq_1~dq_2 \dots  dq_{m-1} ~ dq_m = Q_m(1; (1, -\Bar{Q}_{1, c}, \dots, -\Bar{Q}_{i - 1, c} )) - \Bar{Q}_{m, c},
    \]
which in turn implies that
    \begin{align*}
        \gamma^{\mathcal{I}}(x) &= 1 - m! (Q_m(1; (1, -\Bar{Q}_{1, c}, \dots, -\Bar{Q}_{i - 1, c} )) - \Bar{Q}_{m, c}) \\
    &= 1 - m! \left(
    \frac{1}{m!}(1-x_m^m) -
    \sum_{i=1}^{m-1} \Bar{Q}_{i, c} \frac{1}{(m - i)!} (1 - x_m^{m-i})
    \right) \\
    &= x_m^m + \sum_{i=1}^{m-1} \Bar{Q}_{i, c} \frac{m!}{(m - i)!} (1 - x_m^{m-i}).
    \end{align*}
Now, assume that $c < m$. Then the expression is identical to the one given in Equation \eqref{eq:integral}, with the exception that the lower bound on all integrals after the $c$'th from the inside out will be $x_c$.
The CDF is therefore
    \[\label{eq:rtTMTI}
    \gamma_{\infty, c}(x) = x_c^m + \sum_{i=1}^{m - 1} \bar{Q}_{i, c} \frac{m!}{(m - i)!} (1 - x_c^{m-i}).
    \]
Now, let $c$ be a random variable given by $c = \max\{i \in \{1, \dots, m\}: p_{(i)} < \tau \}$. Assume without loss of generality that $p_1 \leq \dots \leq p_m$. For any fixed integer $i \leq m$ we note that $p_1, \dots, p_i \mid c = i \sim U(0, \tau)$ and that the joint distribution of $(p_{(1)}, \dots, p_{(i)})$ conditional on $c = i$ has density $i! / \tau^i$ on the simplex $\{\boldsymbol{p} \in [0, 1]^i: p_1 \leq \dots \leq p_i < \tau \}$. By the Law of Total Probability, we can write
	\[
	\gamma_{\infty, c} (x) = \sum_{i = 0}^m \prop(Z_{\infty, c} < x \mid c = i) \prop(c = i).
	\]
The distribution of $c$ is binomial with probability parameter $\tau$ and size $m$, i.e., $$\prop(c = i) = \binom{m}{i} \tau^i (1 - \tau)^{m - i} .$$ 
Consider first the case of $i \geq 1$.
We see that
	\begin{align*}
	\prop(Z_{\infty, i} \leq x \mid c = i) &= 1 - \prop(P_1 > x_1, \dots,  P_i > x_i \mid c = i) \\
	& = 1- \frac{i!}{\tau^i} \int_{x_i}^{\tau} \int_{x_{i-1}}^{q_i} \dotsc \int_{x_{2}}^{q_3} \int_{x_1}^{q_2}~
    dq_1~dq_2 \dotsc  dq_{i-1} ~ dq_i I(x_i < \tau) \\
	& = 1 - \frac{i!}{\tau^i} \left( 
		\tilde{Q}_i - \bar{Q}_{i, m}
	\right) I (x_i < \tau).
	\end{align*}
In the case of $i = 0$, we see that $\prop(Z_{\infty, c} < x \mid c = 0) = \prop(\beta(1, m)(P_1) < x \mid c = 0)$. The distribution of $P_1$ conditionally on no $p$-values falling below $\tau$ is uniform on the interval $(\tau, 1)$. Thus, we have $\beta(1, m)(P_1) \mid c = 0 \sim U(\beta(1, m)(\tau), 1)$ and therefore
	\[
	\prop(Z_{\infty, c} < x \mid c = 0) = \frac{x - \beta(1, m)(\tau)}{1 - \beta(1, m)(\tau)} I(x_1 > x).
	\]
	Combining all of the above, we obtain
	\[
	\begin{split}
	\gamma_{\infty, c}(x) &= 
	(1 - \tau)^m \frac{x - \beta(1, m)(\tau)}{1 - \beta(1, m)(\tau)} I(x_1 > x) \\
	& + 
		\sum_{i = 1}^{m} \left[
			\binom{m}{i} \tau^i (1 - \tau)^{m - i} \left\{
				1 - \frac{i!}{\tau^i} \left(
					\tilde{Q}_i
					-
					\bar{Q}_{i,m}
				\right) 
				I(x_i \leq \tau)
			\right\}
		\right]
	\end{split}
	\]
	which proves the claim.

\subsection{Proof of Lemma \ref{lemma:minSubset}}
	Assume without loss of generality that the $p$-values $\boldsymbol{p}_\mathcal{I}$ are already sorted, i.e., that $p_{1} \leq \dots \leq p_{m}$. Then $\eta$ is the identity function and can therefore be omitted entirely. 
	 Fix a set $\mathcal J \coloneqq \{j_1, \dots, j_k\} \in \mathcal{J}^k$ and assume without loss of generality that $j_1 < \dots < j_k$. Fix $j \in \mathcal J$ and $l \in \mathcal I \backslash \mathcal J$ such that $l < j$ and let $\mathcal J_l^{-j} \coloneqq (\mathcal J \backslash \{j\}) \cup \{l\}$, i.e., $\mathcal J_l^{-j}$ is the set obtained by substituting $j$ for $l$ in $\mathcal J$. It suffices to show that
	\begin{displaymath}
		h(\boldsymbol{Y}^{\mathcal{J}_l^{-j}}) \leq h(\boldsymbol{Y}^{\mathcal{J}}).
	\end{displaymath}
	There are two cases to consider:
	
	\textbf{Case 1; $j = \min \mathcal J$:} If $j$ is the smallest element in $\mathcal J$ then substituting it for $l$ does not change the ordering of $\mathcal J$. By Condition \eqref{eq:assumption_weak}, it holds that
	\begin{displaymath}
		F_{(1)}(p_{(l)}) =: Y_{1}^{\mathcal{J}_l^{-j}} \leq Y_{1}^{\mathcal{J}} \coloneqq F_{(1)}(p_{(j)}).
	\end{displaymath}
	As all other values for $\boldsymbol{Y}^{\mathcal{J}_l^{-j}}$ and $\boldsymbol{Y}^{\mathcal{J}}$ are unchanged, it follows from Condition \ref{eq:assumption_h} that $h(\boldsymbol{Y}^{\mathcal{J}_l^{-j}}) \leq h(\boldsymbol{Y}^{\mathcal{J}})$.
	
	\textbf{Case 2; $j > \min \mathcal J$:} Define by $\tilde{j}$ the largest index in $\mathcal J$ smaller than $j$, i.e.
	\[
	\tilde{j} \coloneqq \max\{i \in \mathcal J: i < j\}.
	\]
	Suppose first that $\tilde{j} < l < j$. In this case, the ordering of $\mathcal J$ is unchanged when substituting $j$ for $l$, making this case isomorphic to case 1. If, on the other hand, $l < \tilde {j}$ the ordering of $\mathcal J$ changes when substituting $j$ for $l$. Let $\hat{j}$ be the smallest index in $\mathcal J$ larger than $l$, i.e.,
	\[
		\hat{j} \coloneqq \min \{
			i \in \mathcal J: i > l
		\}.
	\]
	Then we must show two things
	\begin{enumerate}
		\item $Y_{l}^{\mathcal{J}_l^{-j} } \leq Y_{\hat{j}}^{\mathcal{J}}$.
		\item For all $i > \hat{j}$ it holds that $Y_i^{\mathcal{J}_l^{-j}} \leq Y_i^{\mathcal J}$.
	\end{enumerate}
	Let $\eta_\mathcal{A}: \mathcal A \to \mathcal A$ be a function sorting the elements of a set $\mathcal A$. That is, $\eta_\mathcal{A}(a)=b$ if and only if $a$ is the $b$'th lowest element in $\mathcal A$. We then see that $\eta_{\mathcal{J}_l^{-j}}(l) = \eta_{\mathcal{J}}(\hat{j})$ and thus that the first point above is satisfied as we have $p_{(l)} \leq p_{(\hat{j})}$ and therefore $Y_{l}^{\mathcal{J}_l^{-j} } \leq Y_{\hat{j}}^{\mathcal{J}}$ by Condition \eqref{eq:assumption_weak}. The second point above is satisfied by Condition \eqref{eq:assumption_F}, as $\eta_{\mathcal{J}_l^{-j}}(h) = \eta_\mathcal{J}(h) + 1$ for any $h > \hat{j}$. This proves Lemma \ref{lemma:minSubset}. 

\subsection{Proof of Theorem \ref{theorem:main0}}
	We start by reminding the reader that all three TMTI statistics have $F_{(i)}(x) = \beta(i, m + 1 - i)(x)$, for $i = 1, \dots, m$, regardless of the choice of $n$. These functions are weakly increasing, as they are CDFs, thus satisfying Condition \eqref{eq:assumption_weak}. Next, fix $x \in (0,1)$ and $i \in \mathcal{I}$ with $i<m$. We then see that
	\begin{align*}
		F_{(i+1)}(x) &= \sum_{h = i + 1}^k \binom{k}{h} x^h (1-x)^{k-h} \\
 		&=  \beta(i + 1, k + 1 - (i + 1))(x) \\
		&< \binom{k}{i} x^i (1-x)^{k-i} + \sum_{h = i + 1}^k \binom{k}{h} x^h (1-x)^{k-h} \\
 		& = \sum_{h = i}^k \binom{k}{h} x^h (1-x)^{k-h} \\
 		& = \beta(i, k + 1 - i)(x) \\
		& = F_{(i)}(x),
	\end{align*}
	by using Equation \eqref{eq:pbeta_alt}. Thus, Condition \eqref{eq:assumption_F} is satisfied.
	
	Let $\mu > 1$ and note that
	\begin{align*}
		h_{TMTI_\infty}(\boldsymbol{Y})  &\coloneqq \min \boldsymbol{Y} \\
		h_{tTMTI_\infty}(\boldsymbol{Y}) &\coloneqq \min (Y_1 + \mu I(\beta^{-1}(1,m)(Y_1) > \alpha), \dots, Y_m + \mu I(\beta^{-1}(m,1)(Y_m) \geq \alpha)) \\
		h_{rtTMTI_\infty}(\boldsymbol{Y}) &\coloneqq \min (Y_1, \dots, Y_K).
	\end{align*}
	It is immediate that both $h_{TMTI_\infty}$ and $h_{rtTMTI_\infty}$ satisfy Condition \eqref{eq:assumption_h}, as the mapping $\boldsymbol{x} \mapsto \min \boldsymbol{x}$ is weakly increasing in every coordinate. To see that $h_{tTMTI_\infty}$ satisfies Condition \eqref{eq:assumption_h}, note that for fixed $i \in \mathcal I$ it holds that $\beta^{-1}(i, m + 1 - i)$ is strictly increasing, thus making $x \mapsto x + \mu I(\beta^{-1}(i, m + 1 - i)(x) > \alpha)$ a weakly increasing mapping and therefore also making $h_{tTMTI_\infty}$ a weakly increasing mapping.

\subsection{Proof of Theorem \ref{theorem:main}}
Assume without loss of generality that $\mathcal{I} = \{1, \dots, m \}$ and $p_1 \leq \dots \leq p_m$ and denote by $\mathbb{P}_{h(\boldsymbol{Y}^\mathcal{J})}$ the CDF of $h(\boldsymbol{Y}^\mathcal{J})$. The adjusted $p$-value for the test of some $H_{i}$ is the maximal $p$-value across all intersection hypotheses in the closure of $\{i\}$ in $\mathcal{I}$, i.e.
    \[
    p_i^\star \coloneqq \max_{\mathcal{J} \in \{i\}_{\mathcal{I}}^*} \mathbb{P}_{h(\boldsymbol{Y}^\mathcal{J})} \circ h(\boldsymbol{Y}^\mathcal{J}).
    \]
Let $\mathcal{J}_{\{i\}_\mathcal{I}^*}^k$ denote the set of all sets in $\{i\}_\mathcal{I}^*$ of size $k$. Then
    \[
    \argmax_{\mathcal{J} \in \mathcal{J}_{\{i\}_\mathcal{I}^*}^k} \mathbb{P}_{h(\boldsymbol{Y}^\mathcal{J})} \circ h(\boldsymbol{Y}^\mathcal{J}) =
    \begin{cases}
        \{m - k, \dots, m\}, & \text{if } i \geq m - k\\
        \{i, m - k + 1, \dots, m\}, & \text{else}
    \end{cases}
    \]
by Lemma \ref{lemma:minSubset}. Let $\overline{\{i\}_{\mathcal{I}}^*} \coloneqq \{\{i, m\}, \{i, m-1, m\}, \dots, \{i, \dots, m\}, \{i-1, \dots, m\}, \dots, \{1,\dots, m\} \}$. Then, by Lemma \ref{lemma:minSubset}
    \[
    p_i^\star = \max_{\mathcal{J} \in \overline{\{i\}_{\mathcal{I}}^*}} \mathbb{P}_{h(\boldsymbol{Y}^\mathcal{J})} \circ h(\boldsymbol{Y}^\mathcal{J}).
    \]
Thus, the adjusted $p$-value for any hypothesis $H_i$ can be obtained in $\vert\overline{\{i\}_\mathcal{I}^*}\vert = m - 1$ steps\footnote{Disregarding the first step, as a marginal $p$-value for the test of $H_i$ is already supplied.}. However, note for any $i$ that $\vert\overline{\{1\}_\mathcal{I}^*} \cap \overline{\{i\}_\mathcal{I}^*} \vert = i - 1$. Therefore, the number of steps required to obtain an adjusted $p$-value for all hypotheses is $\sum_{i = 1}^{m} (m - i) = \frac{1}{2}m (m -1)$, as claimed.  

\subsection{Proof of Lemma \ref{lemma:skipM}}
Assume without loss of generality that $\mathcal{J}_1 = \{1, \dots, m_1\}$ and $\mathcal{J}_2 = \{1, \dots, m_2\}$. Then $m_1 < m_2$. Let $L \coloneqq \argmin_{j \in \mathcal{J}_1} Y_j^{\mathcal{J}_1}$. Then
    \[
    Z^{\mathcal{J}_1} = \beta(L, m_1 + 1 - L)(P_{(L)}) < \beta(L, m_2 + 1 - L)(P_{(L)}),
    \]
which implies that $\min \boldsymbol{Y}^{\mathcal{J}_2} < \min \boldsymbol{Y}^{\mathcal{J}_1}$ and therefore $\prop(\min \boldsymbol{Y}^{\mathcal{J}_1} < x) < \prop(\min \boldsymbol{Y}^{\mathcal{J}_2} < x)$, which is the claimed result.

\clearpage
\section{Further simulation studies}
\subsection{An investigation of the robustness of the TMTI CDFs against different dependency structures}\label{appendix:dependency_simulation}
In this section we investigate how well the analytical expression of the TMTI CDFs under the null distribution derived in \cref{sec:analytic_computation} under an i.i.d.\ assumption approximates the actual CDF of the TMTI statistics under different dependency structures. 

In the following, we let $\boldsymbol{X} \sim \mathcal{N}(\boldsymbol{0}, \Sigma)$ be an $m$-dimensional random vector with coordinate means of zero, where $\operatorname{diag}(\Sigma) = (1~\dots~1)$. We then calculate $p$-values $\boldsymbol{P}$ as
	\[
	\forall i \in \{1, \dots, m\}: \hspace{1cm}
	P_i \coloneqq 2 (1 - \Phi(\vert X_i \vert)),
	\]
where $x \mapsto \Phi(x)$ is the CDF of a $N(0, 1)$ distribution. This ensures that each $P_i$ is uniform on $(0, 1)$ and that the dependency structure of $\boldsymbol{P}$ is fully determined by the covariance matrix $\Sigma$. We consider three structures of $\Sigma$:
\begin{enumerate}
	\item Equicorrelated tests, where $\Sigma_{i,j} = \rho I(i \neq j) + I(i = j)$, for all $i, j$ and some $\rho \in (0, 1)$.
	\item Block-diagonal tests, where
	\[
		\Sigma = \begin{pmatrix}
		\Sigma_1 & \dots & 0 \\
		\vdots & \ddots & \vdots \\
		0 & \dots & \Sigma_g
		\end{pmatrix}
	\]
	for some $g < m$ and $\Sigma_1, \dots, \Sigma_g$ are themselves equicorrelated with parameter $\rho$.
	\item Autoregressive tests, where $\Sigma_{i, j} = \rho^{\vert i - j \vert}$.
\end{enumerate}
The first point above can happen in a scenario like the simulation study performed in \ref{sec:power} where we combine $T$-tests performed on independent variables but where the standard error is estimated on the basis of a number of covariates. 
It is unlikely to that the correlation between the tests is high, but it is none-zero. Nevertheless, we try this scenario even for large values of $\rho$ in order to investigate what happens under extreme dependencies. 

The second point in the above represents a scenario in which we have performed multiple tests within $g$ groups or individuals that are independent from one another, but where tests performed within the same group are not independent. The last point represents, for example, a design where the tests are spatially correlated, where dependence is highest for neighboring plots.

\begin{figure}[H]
  \centering
  \includegraphics[width = .9\linewidth]{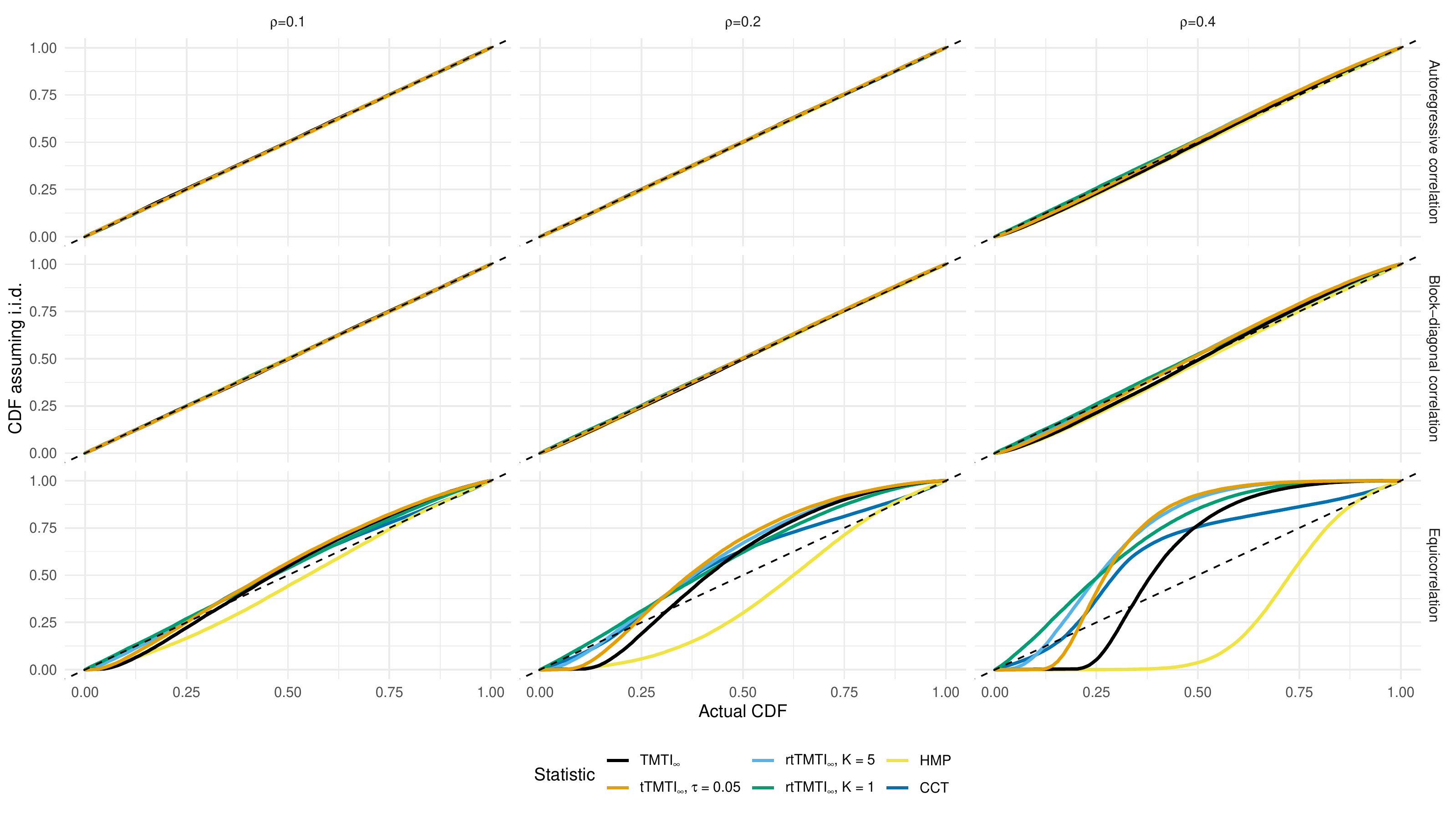}
  \caption{Calibration curves for the CDFs of different TMTI statistics under an i.i.d.\ assumption versus different dependency structures. CCT = {\CCT},  HMP = {\HMP}.}
  \label{fig:dependencyplot}
\end{figure}

\begin{figure}[H]
  \centering
  \includegraphics[width = .9\linewidth]{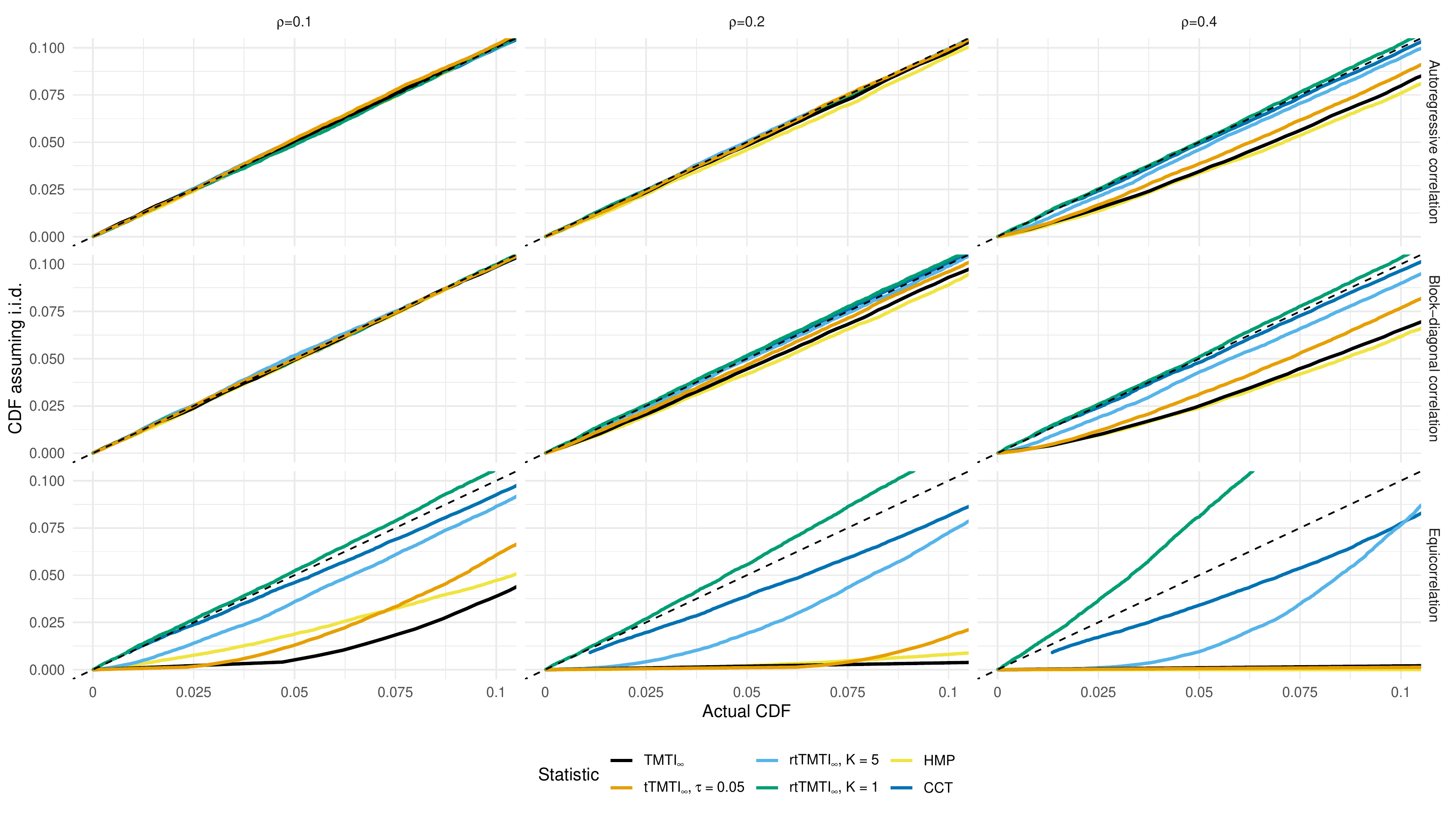}
  \caption{Calibration curves for the CDFs of different TMTI statistics under an i.i.d.\ assumption versus different dependency structures, zoomed to only show the region $(0, 0.1) \times (0, 0.1)$. CCT = {\CCT},  HMP = {\HMP}.}
  \label{fig:dependencyplot_zoom}
\end{figure}

We perform the experiments with $m = 200$ tests and different values of $\rho$. For the block-diagonal experiment we set $g = 40$, corresponding to $40$ groups of five tests each. In each experiment we bootstrap the CDFs of TMTI$_\infty$, tTMTI$\infty$ with $\tau = 0.05$ and rtTMTI$_\infty$ with $K = 5$ and $K = 1$ both under an assumption of i.i.d.\ tests and under the actual dependency structure. For comparison, we have also included the {\CCT} and {\HMP}  tests. We then plot calibration curves, i.e., the curve $x \mapsto (\text{actual } \operatorname{CDF}(x), \text{i.i.d.\ } \operatorname{CDF}(x))$. If the i.i.d.\ CDF is robust against departures from independence then this curve will lie exactly on the diagonal of the unit square. If the i.i.d.\ CDF is conservative it will lie above the diagonal of the unit square and if it is anti-conservative it will lie below the diagonal of the unit square. The results are presented in \cref{fig:dependencyplot} and again in \cref{fig:dependencyplot_zoom} where we have zoomed in on the square $(0, 0.1) \times (0, 0.1)$. From the figures, we see that weak dependencies generally do not affect the CDFs of the TMTI statistics by much but stronger dependencies have a large effect on the CDFs. This is similar to what we see for the {\CCT} and {\HMP}  statistics. Although both of these are claimed to be robust against dependence \citep{acat_original,wilson2019harmonic}, we see that both of these also become anti-conservative in their lower tails when there is sufficiently strong dependence. The {\HMP} appears to generally be more anti-conservative than the TMTI statistics, while the {\CCT} performs slightly better than rtTMTI$_\infty$ with $K = 5$, but worse than rtTMTI$_\infty$ with $K = 1$. In particular, we see that when dependencies are strong, TMTI$_\infty$ is very anti-conservative in its lower tail implying a loss of Type I error control. However, its truncated variants, tTMTI$_\infty$ and rtTMTI$_\infty$, are less affected by the dependencies and are only slightly anti-conservative in many cases with strong dependencies. In particular, rtTMTI$_\infty$ with $K = 1$, corresponding to only using the smallest $p$-value, is very robust against most dependencies and, in contrast to the other CDFs presented, never anti-conservative. Only in the equicorrelated experiment when correlations are strong does it deviate from the unit square diagonal and here it is conservative. We also note that the CDFs of all statistics are conservative at their upper tails. However, it is generally only interesting to see how conservative or anti-conservative the i.i.d.\ approximations are in the regions around typical significance level as conservativeness or anti-conservativeness here can be the difference between making a Type I or II error and making no error.

These experiments suggest that rtTMTI$_\infty$ with $K = 1$ is robust against departures from independence and is thus reasonable to use when making no assumptions on the dependency structure of the $p$-values in question. 
In contrast, the other TMTI statistics can be used under some dependency structures if we believe that the dependencies are not too strong, and the truncated TMTI$_\infty$ variants are generally more robust against dependency than TMTI$_\infty$.

\changed{%

\subsection{An example of applying the TMTI to non-independent data}\label{sec:power}
In this section, we give an example of how the TMTI can be applied to non-independent data in a particular setting.
We consider a variation of the simulation study conducted in \cref{sec:newpower}, where the marginal $p$-values now come from dependent $T$-scores instead of independent $Z$-scores.
Furthermore, we investigate the power of this procedure under different alternatives by means of simulation. Throughout this section, we consider a significance level of $\alpha = 0.05$.

We consider the following setup:
let $\boldsymbol{X}$ be an $n\times g$-dimensional binary matrix of full rank satisfying that every row-sum of $\boldsymbol{X}$ equals one and every column-sum of $\boldsymbol X$ equals $k$, for some $k \in \mathbb N$. Fix a vector $\mu \in \mathbb R^g$ and let $\epsilon \sim N_n(0, I_n)$ be a random variable, where $I_n$ is the $n\times n$ identity matrix . We then define
\[
	\boldsymbol W \coloneqq \boldsymbol X\mu + \epsilon.
\]
The interpretation of this experiment is that we have recorded $k$ observations of some random variable, $W$, in $m$ different groups to obtain a total of $N=mk$ samples. By altering the number of non-zero elements in $\mu$ we control the number of groups that affect the outcome, $\boldsymbol W$. As the power of any test in this scenario is directly associated with the magnitude of the coefficients, $\mu$, we will only consider constant values of $\mu$. That is, the elements in $\mu$ which we allow to be non-zero, will all be equal. We then consider the global hypothesis
$
	H_0: \bigcap_{i = 1}^m (\mu_i = 0),
$
which is the hypothesis that the outcome $W$ is not affected by any of the $m$ groups. The goal is now to estimate the power of the TMTI, i.e., the quantity
$
\prop_{H_A}(\text{reject } H_0),
$ 
under different alternative hypotheses $H_A$, i.e., different $\mu$ vectors. For every alternative hypothesis considered here, we will compute the $p$-value under 
four tests
from the TMTI family; 
TMTI$_\infty$ and its truncated and rank truncated versions (with $\tau = 0.05$ and $K = 5$), and TMTI$_1$. For TMTI$_1$, we do not consider any truncated variants, as we expect these to be roughly equal (per \cref{fig:gamma_comparison}).

For comparison, we consider a number of combination-based test procedures. For this, we will compute the marginal $T$-tests of the hypotheses $H_i: \mu_i = 0$ for $i = 1, \dots, m$ under the joint model. Then, we will apply the following procedures; the {\CCT}  \citep{acat_original}; the {\HMP}  \citep{wilson2019harmonic}; the {\CBAM}  \citep{vovk2020combining}; a Bonferroni correction; and the {\rTPM}.
Additionally, we include the $F$-test, as this would often be the natural choice of test in this particular setup. However, the $F$-test offers less flexibility compared to a combination test, e.g., when used in a Closed Testing Procedure (see \cref{sec:MTP}), as this would require fitting $2^m$ linear models.

In order to calculate $p$-values for the TMTI statistics, we estimate the $\gamma$ functions under $H_0$ by employing a bootstrapping scheme: For $i =1,\dots, 10^5$ we first define $\tilde{\boldsymbol{W}} \coloneqq \boldsymbol{W} - (X^TX)^{-1}X^T \boldsymbol{W} X^T$. That is, $\tilde{\boldsymbol W}$ is $\boldsymbol{W}$ with the group means subtracted. This ensures that each group of $\tilde{\boldsymbol W}$ has mean zero. We then construct $\boldsymbol{W}_i$ by resampling $\boldsymbol{W}$ uniformly with replacement in order to introduce variation. We then compute the relevant marginal $T$-test statistics, compute $p$-values and finally output the TMTI statistic. We then use the empirical CDFs of the bootstrapped TMTI statistics.
In \cref{fig:size_simulation} we show the simulated sizes of all included tests.
From this we conclude that all tests have approximately the correct size or lower, except for the {\FCT} and the {\TPM}  (both of which assume independence of the $p$-values), and these are therefore left out from further simulations. One can apply the same bootstrapping scheme as described above to obtain valid $p$-values for both of these tests. However, we do refrain from doing that here. We note also that the {\rTPM} has a slightly increased Type I error, but it is sufficiently little that this may be attributed to chance. The increased variance in the estimates of the Type I errors for the TMTI statistics is due to the critical values of these being estimated by bootstrapping.

\begin{figure}
    \centering
    \includegraphics[width = .9\linewidth]{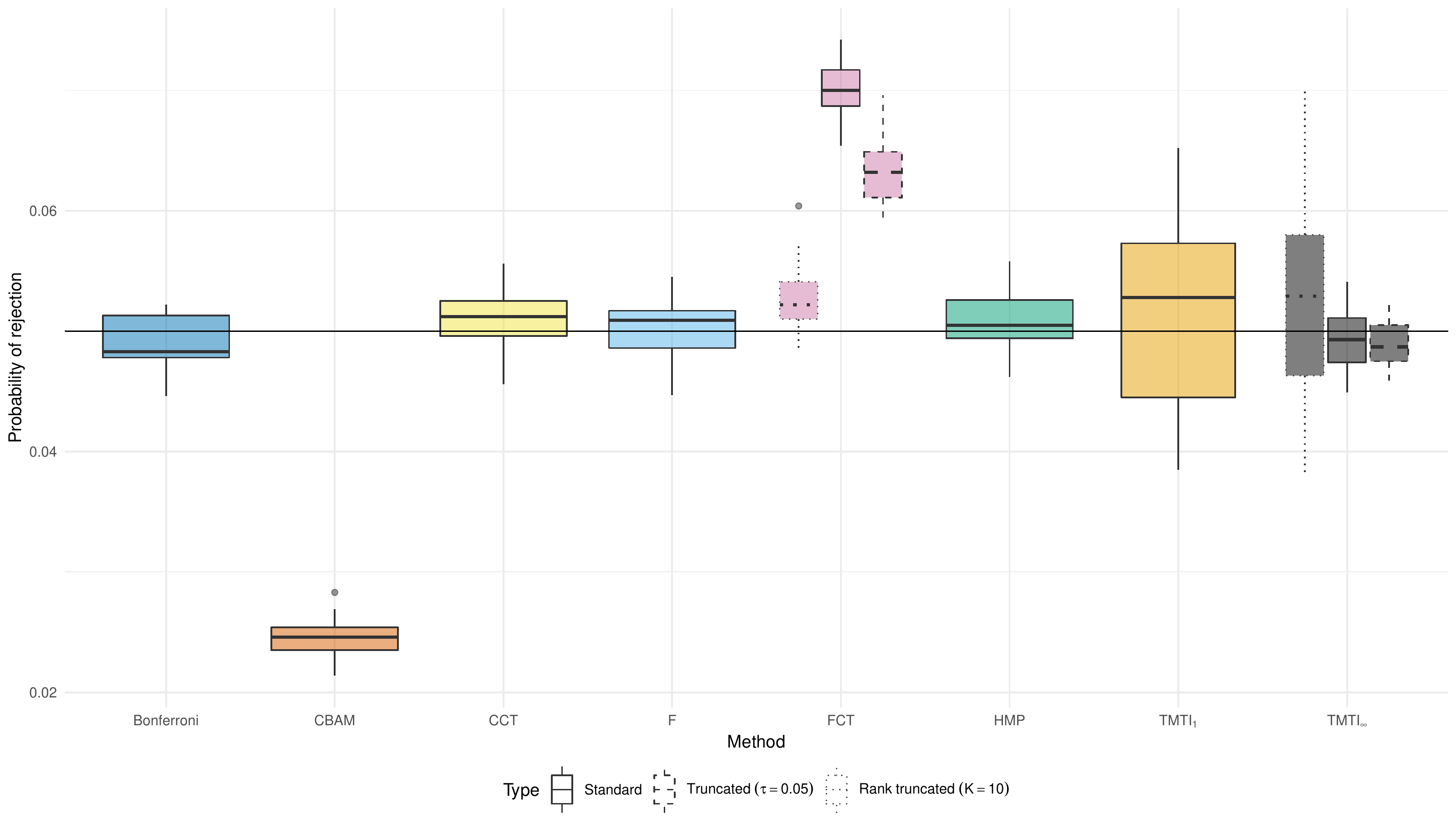}
    \caption{Estimated sizes of the included tests in \ref{sec:power}, i.e., the probability of rejecting the global null hypothesis when it is true. The \TPM and \rTPM are denoted as the truncated FCT and rank truncated FCT, respectively.
    } %
    \label{fig:size_simulation}
\end{figure}

\cref{fig:dependent_power} contains the results of a simulation with $N_{\operatorname{false}} \in \{10^0, \dots, 10^3\}$. Generally, the results are similar to those displayed in \cref{fig:new_simulation}: TMTI$_\infty$ and tTMTI$_\infty$ both perform well in all settings, although not as well as some methods when $N_{\operatorname{false}} = 1$, but still better than the $F$-test. When $N_{\operatorname{false}} > 1$, we generally find that TMTI$_\infty$ and tTMTI$_\infty$
perform either as good or better than the best of the competing methods. TMTI$_1$ and rtTMTI$_\infty$ have similar performance: when there are less than ten false hypotheses, these perform on par or better than with the best of the competing methods, but when there are multiple false hypotheses, they slightly outperform the {\CCT}  and {\HMP},  although they are both outperformed by the $F$-test, TMTI$_\infty$ and tTMTI$_\infty$.

\begin{figure}
    \centering
    \includegraphics[width = .9\linewidth]{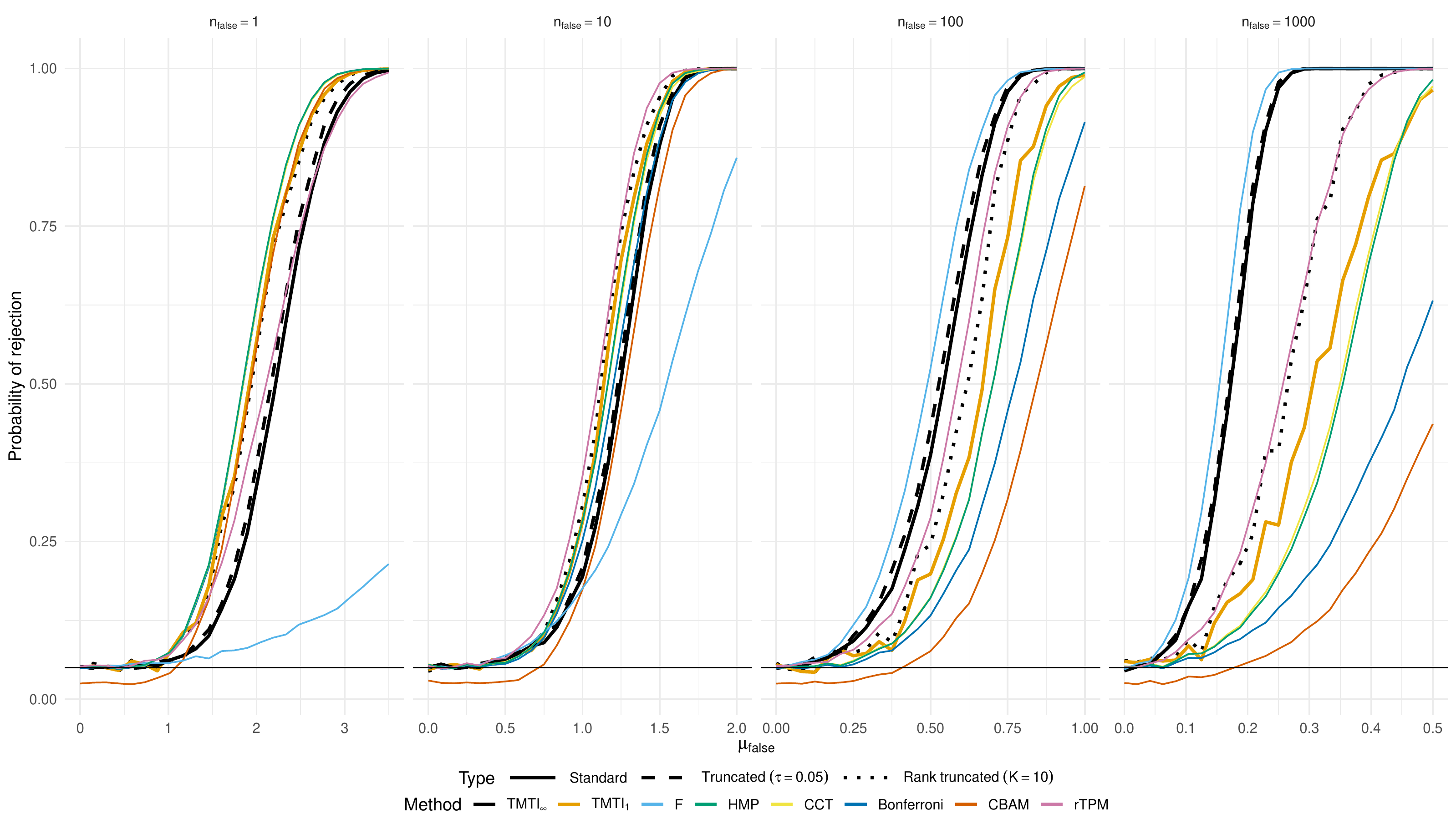}
    \caption{Power curves for different TMTI tests and competing methods. Generally, TMTI$_\infty$ and tTMTI$_\infty$ work well in all settings, while TMTI$_1$ and rtTMTI$_\infty$ performs best when signals are sparse.}
    \label{fig:dependent_power}
\end{figure}

As in \cref{sec:newpower}, the results of these simulations indicate that TMTI$_\infty$ and tTMTI$_\infty$ offer an alternative to current tests that is powerful against a wide range of alternative hypotheses. This is useful if one has no \textit{a priori} knowledge of the sparsity and strength of signals, as well as when employing the tests in a {\CTP}. If one has \textit{a priori}  knowledge, that the signals are sparse and strong, one should rather employ a test such as TMTI$_1$, rtTMTI$_\infty$, {\HMP}, {\CCT} or a Bonferroni test.

\begin{remark}
	The above example generalizes to a situation in which we have two sets of covariates, $\boldsymbol{X}$ and $\tilde{\boldsymbol{X}}$, such that 
    $
		\boldsymbol{W} \coloneqq \boldsymbol{X} \mu + \tilde{\boldsymbol{X}} \tilde{\mu} + \epsilon,
	$
	where $\tilde{\boldsymbol{X}}$ are covariates we would simply like to adjust for when computing our test, but not variables that are of interest. In this setting, we could also compute the $p$-values corresponding to the tests of $H_i: \tilde{\mu}_{i} = 0$, although we would not care whether they are true or not. Thus, we would select a subset of the $p$-values, say $\mathcal J$, corresponding to those related to $\mu$, and compute our test only for those.
\end{remark}
}

\changed{%
\subsection{The effects of mixed $\mu$ values}\label{appendix:mixed_betas}
In this section, we repeat the experiment performed in \cref{sec:newpower}, with the change that the $\mu$ values are allowed to differ between the false marginal hypotheses. Thus, for $N_{\operatorname{false}}$ false hypotheses, we generate $p$-values by sampling $X_{\operatorname{false}, i} \sim N(\mu / i, 1)$ (i.e., the signal with the largest effect has mean $\mu$ and the signal with the weakest effect has mean $\mu / N_{\operatorname{false}}$), where $\mu$ is chosen equidistantly between the values that satisfy that a Bonferroni test has either 5\% or 99\% power to reject the global null hypothesis in a setting with no conservatism. The results are displayed in \cref{fig:mixed_betas}. The comments to this figure are the same as the comments to \cref{fig:new_simulation}.

\begin{figure}
    \centering
    \includegraphics[width = .9\linewidth]{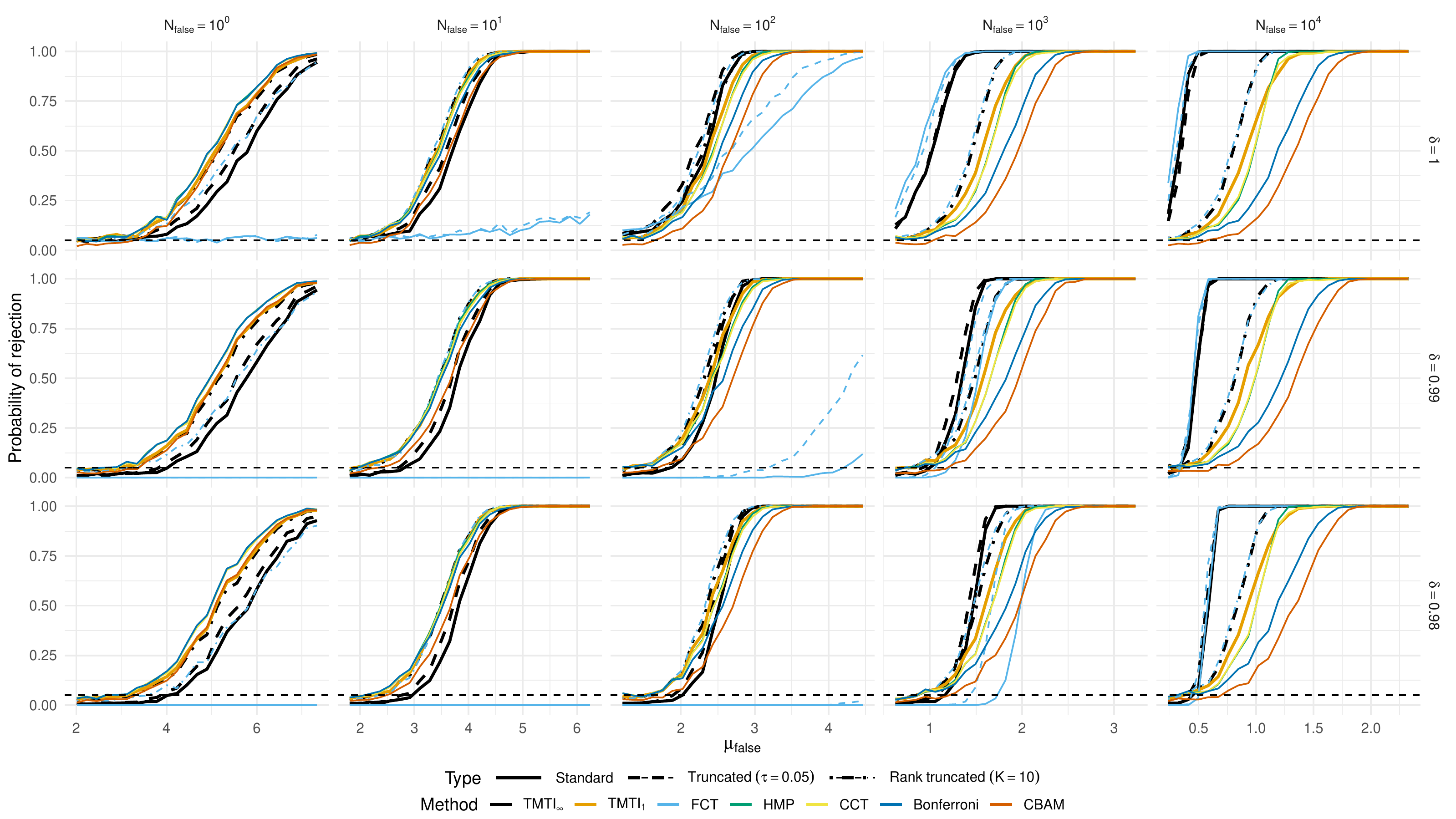}
    \caption{Power curves for different TMTI tests and competing methods, in a scenario where the signal strength is allowed to differ between each false marginal hypothesis. The values of $\mu_{\operatorname{false}}$ are chosen equidistantly between the two values, which satisfy that a Bonferroni test has either 5\% or 99\% power to reject the global null hypothesis in a setting with no conservatism.}
    \label{fig:mixed_betas}
\end{figure}
}

\section{Table of adjusted $p$-values for all tests employed in \cref{sec:realdata}}\label{appendix:pvalue_table}
\begin{table}[t]
\begin{tabular}{p{2.5cm}p{1.75cm}p{1.75cm}p{1.75cm}p{1.75cm}p{1.75cm}p{1.25cm}}
\hline
State (change) & $p$-value & Bonferroni & TMTI$_\infty$& rtTMTI$_\infty$/ TMTI$_\infty$ & FCT& rTPM/ FCT~\\ \hline
GA    (-0.323)                   & 0.85628  &1.00000   &0.87219              &0.93682            &0.85753              &0.93775   \\
AR    (-0.777)                   & 0.60282  &1.00000  &0.87219              &0.93682            &0.85753              &0.93775   \\
AL    (-1.568)                   & 0.44008  &1.00000  &0.85873              &0.93682            &0.81333              &0.93775   \\
NJ    (1.565)                    & 0.41998  &1.00000  &0.85873              &0.93682            &0.80157              &0.93775   \\
NE    (1.334)                    & 0.38640  &1.00000  &0.85873              &0.93682            &0.78021              &0.93775   \\
ND    (1.526)                    & 0.36890  &1.00000  &0.85873              &0.93682            &0.76813              &0.93775   \\
DE    (1.374)                    & 0.31162  &1.00000  &0.85873              &0.92675            &0.72551              &0.92768   \\
MI    (2.215)                    & 0.23522  &1.00000  &0.80175              &0.88412            &0.66845              &0.88500   \\
LA    (2.637)                    & 0.20964  &1.00000  &0.78923              &0.88412            &0.64602              &0.88500   \\
IN    (2.149)                    & 0.19388  &1.00000  &0.78923              &0.88412            &0.63076              &0.88500   \\
WI    (2.801)                    & 0.15872  &1.00000  &0.78923              &0.85060            &0.59172              &0.85144   \\
VA    (2.858)                    & 0.14374  &1.00000  &0.77357              &0.84467            &0.57388              &0.84550   \\
WV    (2.331)                    & 0.10026  &1.00000  &0.68933$^\ddagger$   &0.74677            &0.51177$^\ddagger$   &0.74750   \\
MD    (3.339)                    & 0.08226  &1.00000  &0.68933$^\ddagger$   &0.69934            &0.48059$^\ddagger$   &0.71026   \\
CA    (3.777)                    & 0.07912  &1.00000  &0.68454$^\ddagger$   &0.70957            &0.47464$^\ddagger$   &0.71026   \\
OH    (3.466)                    & 0.06590  &1.00000  &0.62312$^\ddagger$   &0.64033            &0.44713$^\ddagger$   &0.64096   \\
NY    (4.893)                    & 0.05802  &1.00000  &0.58342$^\ddagger$   &0.59203            &0.42838$^\ddagger$   &0.59262   \\
PA    (4.303)                    & 0.05572  &1.00000  &0.58342$^\ddagger$   &0.57683            &0.42250$^\ddagger$   &0.57739   \\
FL    (3.784)                    & 0.05490  &1.00000  &0.58342$^\ddagger$   &0.57129            &0.42036$^\ddagger$   &0.57185   \\
WY    (2.226)                    & 0.04678* &1.00000  &0.58342$^\ddagger$   &0.51259            &0.39755$^\ddagger$   &0.51308   \\
NM    (2.334)                    & 0.04650* &1.00000  &0.58342$^\ddagger$   &0.51043            &0.39671$^\ddagger$   &0.51092   \\
CT    (3.204)                    & 0.04104* &1.00000  &0.55925$^\ddagger$   &0.46666            &0.37939$^\ddagger$   &0.46711   \\
OK    (4.181)                    & 0.02036* &0.69224  &0.42037$^\ddagger$   &0.26549            &0.29050$^\ddagger$   &0.26573   \\
KY    (4.326)                    & 0.00964* &0.32776  &0.28899$^\dagger$    &0.13524$^\ddagger$ &0.21234$^\dagger$    &0.13535$^\ddagger$   \\
AZ    (4.993)                    & 0.00904* &0.30736  &0.27561$^\dagger$    &0.12735$^\ddagger$ &0.20643$^\dagger$    &0.12747$^\ddagger$   \\
ID    (2.956)                    & 0.00748* &0.25432  &0.23899$^\dagger$    &0.10651$^\ddagger$ &0.18974$^\dagger$    &0.10659$^\ddagger$   \\
TX    (5.645)                    & 0.00404* &0.13736  &0.17114$^\dagger$    &0.05892$^\dagger$  &0.14480$^\dagger$    &0.05897$^\dagger$   \\
CO    (4.326)                    & 0.00282* &0.09588  &0.12797$^\dagger$    &0.04148*           &0.12286$^\dagger$    &0.04150*   \\
IA    (4.811)                    & 0.00200* &0.06800  &0.11058$^\dagger$    &0.02958*           &0.10453$^\dagger$    &0.02961*   \\
NH    (4.422)                    & 0.00180* &0.06120  &0.10121$^\dagger$    &0.02666*           &0.09939$^\dagger$    &0.02667*   \\
NC    (7.265)                    & 0.00002* &0.00068* &0.00346*             &0.00346*           &0.00843*             &0.00064*   \\
HI    (5.550)                    & 0.00002* &0.00068* &0.00346*             &0.00346*           &0.00843*             &0.00064*   \\
MN    (6.421)                    & 0.00002* &0.00068* &0.00346*             &0.00346*           &0.00843*             &0.00064*   \\
RI    (5.094)                    & 0.00001* &0.00034* &0.00198*             &0.00198*           &0.00551*             &0.00044*   \\ \hline
\end{tabular}
\caption{State-wise changes in mathematics achievements from 1990 to 1992 and the $p$-values for the corresponding $T$-tests. Values that are significant at $\alpha = 0.05$ are marked with an asterisk. Hypotheses that can be rejected while controlling the $k$-FWER at $k=2$ are marked with a $\dagger$. %
}
\end{table}
\end{document}